\documentclass[final]{amsart}
\usepackage{pdfsync}
\usepackage{amsmath,amssymb}
\usepackage{enumerate}
\usepackage{xspace}
\usepackage{boxedminipage}
\usepackage{algorithmicx}
\usepackage[ruled]{algorithm}
\usepackage{algpseudocode}
\usepackage{listings}
\usepackage{tabularx}
\usepackage{subfigure}
\usepackage{tristan}
\usepackage{fullpage}
\renewcommand{\bi}[2]{\ensuremath{\cA\qp{#1,#2}}}
\renewcommand{\bih}[2]{\ensuremath{\cA_h\qp{#1,#2}}}
\renewcommand{\bid}[2]{\ensuremath{\cA^*\qp{#1,#2}}}

\renewcommand{\vec}[1]{\geovec{#1}}
\renewcommand{\mat}[1]{\geomat{#1}}

\renewcommand{\enorm}[2]{\ensuremath{\Norm{#1}}_{dG,{#2}}}
\newcommand{\opnorm}[2]{\ensuremath{\Norm{#1}}_{\widetilde{dG},{#2}}}
\renewcommand{\eenorm}[2]{\,|\!|\!| {{#1}} |\!|\!|_{dG,{#2}}}

\numberwithin{equation}{section}
\setboolean{showtodo}{true}
\setboolean{shownotes}{true}
\setboolean{showchanges}{false}
\setboolean{usemathrsfs}{false}
\author{
  Tristan Pryer
}
\address{
  Tristan Pryer
  \thanks{
    Department of Mathematics and Statistics, Whiteknights, PO Box 220, Reading RG6 6AX, UK
    {\tt{T.Pryer@reading.ac.uk}}.
}}

\title[A posteriori analysis of inconsistent, nonconforming FEMs]
      {An a posteriori analysis of some
        inconsistent, nonconforming Galerkin methods approximating elliptic
        problems}
\date{\today}
\pdfformat{true}

\begin{document}
\maketitle
\begin{abstract}
  In this work we present an 
  a posteriori analysis for classes of inconsistent, nonconforming
  schemes approximating elliptic problems. We show the estimates
  coincide with existing ones for interior penalty type discontinuous
  Galerkin approximations of the Laplacian and give new estimates for
  inconsistent discontinuous Galerkin approximation schemes of
  elliptic problems under quadrature approximation. We also examine
  the effect of inconsistencies on the a posteriori analysis of
  schemes applied to an unbalanced problem.
\end{abstract}

\section{Introduction}
\label{sec:introduction}

The design of efficient numerical approximations of partial
differential equations (PDEs) is of paramount importance in scientific
computing. To this end a large amount of work has been carried out on
finite element methods and, in particular, on the a posteriori
analysis of such methods. The a posteriori framework allows for the
construction of adaptive methods that, given a fixed number of degrees
of freedom, are able to better approximate nonsmooth solutions of PDEs
than the nonadaptive counterparts \cite{Ainsworth:2000,Verfurth:1996}.

Discontinuous Galerkin (dG) methods form a class of finite element
approximation schemes. These, as the name suggest, are constructed
without enforcing continuity of the discrete solution over the
domain. For second order elliptic problems these are nonconforming
methods. See \cite{ArnoldBrezziCockburnMarini:2001} for an accessible
overview and history of these methods for second order problems. For
higher order problems, for example the (nonlinear) biharmonic problem,
dG methods are a useful alternative to using $\cont{1}$ conforming
elements whose derivation and implementation can become very
complicated
\cite{Baker:1977,GeorgoulisHoustonVirtanen:2011,Pryer:2014}. The
historical origins of the method can be traced back to that of the
transport equation and, as such, the dG methods in the non-selfadjoint
case, for example convection dominated problems, are extremely
powerful \cite{Cockburn:1999}.

In the elliptic case there are two main approaches to the a posteriori
analysis of such nonconforming schemes. The first is based on a
Helmholtz decomposition of the error \cite{BeckerHansboLarson:2003}
and the second based on reconstruction operators
\cite{KarakashianPascal:2003}. The reconstruction approach is
extremely versatile and has been applied not just to elliptic problems
but also to those of hyperbolic type
\cite{GeorgoulisHallMakridakis:2014,GiesselmannMakridakisPryer:2014}. In
both scenarios the main idea is to develop an appropriate
\emph{reconstruction} of the dG object which is itself smooth enough
to apply stability arguments of the underlying PDE. The blanket
assumption in the elliptic case is that, up to mesh dependent terms,
the dG approximation is ``close enough'' to a conforming approximation
which is indeed true for primal methods chosen with large enough
penalty parameters \cite{Di-PietroErn:2012}.

In this work we utilise these reconstruction operators to develop an
abstract a posteriori framework for use in elliptic problems. We are
particularly interested in the effect of inconsistencies in the
estimate. The inconsistencies that we on focus may arise in various
forms, quadrature approximation being the most obvious example. This
has been studied previously in the form of data oscillation in the
convergence of $h$-adaptive algorithms
\cite[c.f.]{Morin:2000}. Another related work is
\cite{NochettoSchmidtSiebertVeeser:2006} where the authors study
quadrature effects of which arise from a nonlinear source
term. Inconsistencies may come from other sources though, for example,
stabilisation introduced in the approximation of convection dominated
problems \cite{Burman:2005} may be viewed in this fashion. Mixed
methods like an operator splitting of a nonconstant coefficient fourth
order problem or discretisation of a nonvariational problem
\cite{LakkisPryer:2011a,DednerPryer:2013} are inherently inconsistent
as are those of virtual element type \cite{BC:2013}. With this range
of inconsistencies in mind we build our a posteriori estimate to
resemble Strang's Lemma, the starting point when conducting a priori
error analysis of inconsistent methods. We show that as long as a
reconstruction operator is available, inconsistencies such as those
mentioned above are a posteriori controllable. Interestingly,
depending on the class of inconsistency, it is not always necessary to
a posteriori compute the reconstruction, in some cases it suffices to
only have existence of such an operator.

By building the a posteriori bounds in this way and using ingredients
from standard a posteriori analysis we can deal with problems as
diverse as unbalanced problems like linear nonvariational and
nonlinear $p$--Laplacian or $p$--biharmonic type, although in this introductory
work we will restrict our attention to linear problems. Specifically, to illustrate the
approach, we initially focus on the comparison of two dG
discretisations of the Laplacian which are standard in the
literature. We move onto the a posteriori effect that quadrature error
has on finite element discretisations of a model
problem. We show that the estimator can be split into components of
which some terms can be computed \emph{exactly} using the same
quadrature approximation in the scheme and others requiring further
approximation. We also examine unbalanced elliptic problems that are
nonvariational. We derive $\sobh{2}$ a posteriori bounds which are
shown to be reliable and efficient up to data oscillation terms. We
also examine further effects of quadrature approximation in this
setting and compare with the variational case.

The paper is set out as follows: In \S\ref{sec:setup} we look at an
abstract model problem and develop an a posteriori Strang type
estimate. In \S\ref{sec:selfadj} we then illustrate how this framework
can be utilised for a model problem and consider inconsistencies
arising from quadrature approximations. In \S\ref{sec:unbalanced} we
study the application of unbalanced elliptic problems, these are
nonvariational in nature and we conduct an a posteriori analysis of
existing schemes proposed for this class of problem. Finally, in
\S\ref{eq:numerics} we summarise extensive numerical experiments
testing the robustness of the proposed estimators.

\section{Abstract problem setup}
\label{sec:setup}

Throughout this section we develop an abstract a posteriori framework
for using notation borrowed from \cite[\S 2]{ErnGuermond:2004}. To
that end let $W$ be a Banach space and $V$ a reflexive Banach space,
equipped with norms $\Norm{\cdot}_W, \Norm{\cdot}_V$ respectively. Let
$\cA: W\times V \to \reals$ be a bilinear form and $l:V \to \reals$ be
a linear form. Consider the abstract problem to find $u\in W$ such
that
\begin{equation}
  \label{eq:abstract-prob}
  \bi{u}{v} = l(v) \Foreach v\in V.
\end{equation}

\begin{Hyp}[The continuous problem is inf-sup stable]
  \label{ass:cont-inf-sup}
  Assume that this problem is stable, in the sense that for some
  $\gamma >0$
  \begin{equation}
    \sup_{v\in V,\Norm{v}_V\leq 1} \bi{u}{v} \geq \gamma \Norm{u}_W.
  \end{equation}
\end{Hyp}

\begin{Pro}[A priori bound]
  Under the conditions of Assumption \ref{ass:cont-inf-sup} we have
  the following a priori bound for the solution of (\ref{eq:abstract-prob}):
  \begin{equation}
    \Norm{u}_W \leq C \Norm{l}_{V^*},
  \end{equation}
  where $V^*$ is the \emph{dual space} of $V$.
\end{Pro}
We let $V_h, W_h$ be
finite dimensional (not necessarily conforming) spaces. Suppose that
$V_h, W_h$ are equipped with norms $\Norm{\cdot}_{V_h},
\Norm{\cdot}_{W_h}$ respectively and $W_h$ a seminorm
$\norm{\cdot}_{W_h}$. Let $W(h) = W + W_h$ and let $W(h)$ be equipped
with the norm $\Norm{\cdot}_{W(h)}$ such that
\begin{gather}
  \label{eq:norm-1}
    \Norm{w_h}_{W(h)} = \Norm{w_h}_{W_h} \Foreach w_h \in W_h
    \\
  \label{eq:norm-2}
    \Norm{w}_{W(h)} = \Norm{w}_W \Foreach w \in W
    \\
  \label{eq:norm-3}
    \norm{w}_{W_h} = 0 \Foreach w\in W
    \\
  \label{eq:norm-4}
    \Norm{v}_{V_h} = \Norm{v}_V \Foreach v \in V.
\end{gather}
We consider the Galerkin method to seek $u_h \in W_h$ such that
\begin{equation}
  \label{eq:galerkin-method}
  \bih{u_h}{v_h} = l_h(v_h) \Foreach v_h\in V_h,
\end{equation}
where $\bih{\cdot}{\cdot} : W_h \times V_h \to \reals$ is an
approximation of $\bi{\cdot}{\cdot}$ and $l_h\qp{\cdot}:V_h \to
\reals$ is an approximation of $l\qp{\cdot}$.

\begin{Hyp}[On the Galerkin method]
  \label{ass:galerkin-bdd}
  We will assume that the Galerkin scheme is bounded in the sense that
  \begin{equation}
    \bih{w_h}{v_h} \leq C_B \Norm{w_h}_{W_h}\Norm{v_h}_{V_h} \Foreach w_h\in W_h, v_h \in V_h.
  \end{equation}
\end{Hyp}

\begin{Rem}
  Note that the is \emph{no} assumption on the convergence
  of the Galerkin scheme, only that the problem is
  stable.
\end{Rem}

\begin{Hyp}[Conforming finite dimensional space]
  Suppose another finite dimensional space, $\fesW$, exists and is a
  conforming approximation of $W$, \ie $\fesW\subset W$ and there
  exists a \emph{reconstruction operator} $E:W_h\to \fesW$. 
\end{Hyp}

\begin{The}[Strang type a posteriori result]
  \label{the:the-aposteriori-strang}
  Let $e:= u - u_h = \qp{u - E(u_h)} + \qp{E(u_h) - u_h} =: e^C + e^N$
  then we have the a posteriori result that for any $v_h \in V_h$
  \begin{equation}
    \label{eq:aposteriori-primal}
    \begin{split}
      \Norm{e}_{W_h}
      &\leq
      \frac{1}{\gamma}
      \bigg(
      \sup_{v\in V, \Norm{v}_V \leq 1}
      \qb{
        l_h(v-v_h) - \bih{u_h}{v - v_h}
      }
      +
      \sup_{v\in V, \Norm{v}_V \leq 1}
      \qb{
        l(v) - l_h(v)
      }
      \\ 
      &\qquad \qquad +
      \sup_{v\in V, \Norm{v}_V \leq 1}
      \qb{
        \bih{E(u_h)}{v} - \bi{E(u_h)}{v}
      }
      \bigg)
      +
      \qp{1+ \frac{C_B}{\gamma}} \Norm{e^N}_{W_h}
      \\
      & =: 
      \frac{1}{\gamma} \qp{ \E_1 + \E_2 + \E_3 + \qp{\gamma + C_B }\E_4}.
    \end{split}
  \end{equation}
\end{The}

\begin{Rem}[Structure of the estimate]
  It should be noted that the result in Theorem
  \ref{the:the-aposteriori-strang} is not a posteriori computable
  since it involves taking the supremum over an infinite dimensional
  space. It is, however, the beginnings of an a posteriori bound. The
  first term $\E_1$ represents a standard residual term which can be
  computably estimated using standard techniques
  \cite[c.f.]{Ainsworth:2000}. The second and third terms $\E_2, \E_3$
  measure a form of inconsistency error and the fourth term, $\E_4$,
  the nonconformity error.
\end{Rem}

\begin{Rem}[Comparison with the a priori Strang Lemma]
  Recall the Strang a priori result that if the bilinear form
  associated to the Galerkin method (\ref{eq:galerkin-method})
  satisfies a discrete inf--sup condition, that is
  \begin{equation}
    \sup_{v_h\in V_h,\Norm{v_h}_{V_h} \leq 1} 
    {\bih{w_h}{v_h}}
    \geq 
    \gamma_h \Norm{w_h}_{W_h}
  \end{equation}
  then
  \begin{equation}
    \label{eq:apriori-strang}
    \begin{split}
    \Norm{e}_{W_h} 
    &\leq
    \qp{1 + \frac{C_B}{\gamma_h}}
    \inf_{w_h\in V_h}\Norm{u - w_h}_{W_h}
    +
    \frac{1}{\gamma_h}
    \bigg(
      \sup_{v_h \in V_h, \Norm{v_h}_{V_h} \leq 1}
      \qb{l_h\qp{v_h} - \bih{u}{v_h}}
      \bigg).
    \end{split}
  \end{equation}
  The inconsistency term represents how badly the solution $u$
  satisfies the discrete scheme. Theorem
  \ref{the:the-aposteriori-strang} may be reformulated as
  \begin{equation}
    \begin{split}
      \Norm{e}_{W_h}
      &\leq
      \frac{1}{\gamma}
      \bigg(
      \sup_{v\in V, \Norm{v}_V \leq 1}
      \qb{
        l_h(v-v_h) - \bih{u_h}{v - v_h}
      }
      +
      \sup_{v\in V, \Norm{v}_V \leq 1}
      \qb{
        l(v)  - \bi{E(u_h)}{v}
      }
      \\ 
      &\qquad \qquad +
      \sup_{v\in V, \Norm{v}_V \leq 1}
      \qb{
        l_h(v) - \bih{E(u_h)}{v}
      }
      \bigg)
      +
      \qp{1+ \frac{C_B}{\gamma}} \Norm{e^N}_{W_h}
    \end{split}
  \end{equation}
  and we may think of the inconsistencies as how badly the
  reconstruction $E(u_h)$ satifies both the smooth problem and
  the Galerkin approximation.
\end{Rem}

\begin{Proof}{of Theorem \ref{the:the-aposteriori-strang}.}
  Using the stability of the continuous problem from Assumption
  \ref{ass:cont-inf-sup} we have for any $v_h\in V_h$
  \begin{equation}
    \begin{split}
    \Norm{e^C}_{W} 
    &\leq 
    \frac{1}{\gamma}
    \qp{
    \sup_{v\in V, \Norm{v}_V \leq 1} 
    \bi{e^C}{v} 
  }
    =
    \frac{1}{\gamma}
    \qp{
      \sup_{v\in V, \Norm{v}_V \leq 1} 
    \qb{
      \bi{u}{v} 
      -
      \bi{E(u_h)}{v}
    }}
    \\
    &
    \leq
    \frac{1}{\gamma}
    \qp{
    \sup_{v\in V, \Norm{v}_V \leq 1} 
    \qb{
      l(v) - \bi{E(u_h)}{v}
    }
  }
    ,
  \end{split}
  \end{equation}
  using the continuous problem (\ref{eq:abstract-prob}). Adding and
  subtracting appropriate terms, we have
  \begin{equation}
    \begin{split}
    \Norm{e^C}_{W} 
    &\leq 
    \frac{1}{\gamma}
    \bigg(
    \sup_{v\in V, \Norm{v}_V \leq 1} 
    \bigg[
    l(v) - l_h(v) + l_h(v - v_h) - \bih{u_h}{v - v_h}
      \\ 
      &\qquad\qquad\qquad\qquad  - \bih{E(u_h) - u_h}{v}  + \bih{E(u_h)}{v} - \bi{E(u_h)}{v}
    \bigg]
    \bigg)
    \\
    &\leq 
    \frac{1}{\gamma}
    \bigg(
    \sup_{v\in V, \Norm{v}_V \leq 1} 
    \qb{
      l_h(v - v_h) - \bih{u_h}{v - v_h}
    }
    +
    \sup_{v\in V, \Norm{v}_V \leq 1} 
    \qb{
      l(v) - l_h(v) 
    }
    \\
    &\qquad 
    +
    \sup_{v\in V, \Norm{v}_V \leq 1} 
    \qb{
      \bih{E(u_h)}{v} - \bi{E(u_h)}{v}
    }
    +
    \sup_{v\in V, \Norm{v}_V \leq 1} 
    \bih{E(u_h) - u_h}{v}\bigg),
    \end{split}
  \end{equation}
  in view of the triangle inequality. To conclude we note that
  \begin{equation}
    \sup_{v\in V, \Norm{v}_V \leq 1} 
    \bih{E(u_h) - u_h}{v}
    \leq
    C_B \Norm{e^N}_{W_h} \Norm{v}_{V_h},
  \end{equation}
  from Assumption \ref{ass:galerkin-bdd}, (\ref{eq:norm-2}) and
  (\ref{eq:norm-4}). The result then follows from 
  \begin{equation}
    \Norm{e}_{W_h} \leq \Norm{e^C}_W + \Norm{e^N}_{W_h}, 
  \end{equation}
  concluding the proof.
\end{Proof}

\begin{Rem}[Computability of (\ref{eq:aposteriori-primal})]
  \label{rem:computability}
  For (\ref{eq:aposteriori-primal}) to be useful some assumptions need
  to be made on the reconstruction $E$. In particular, we impose that 
  \begin{equation}
    \Norm{E(w_h) - w_h}_{W_h} \leq C \norm{w_h}_{W_h} \Foreach w_h \in W_h.
  \end{equation}
  Examples of reconstruction operators satisfying this condition for
  various test problems will be given in subsequent chapters.
\end{Rem}

\begin{Defn}[Formal dual problem]
  Let $L$ be a Hilbert space such that for any $v\in V$ we have
  \begin{equation}
    \Norm{v}_L \leq C\Norm{v}_V.
  \end{equation}
  We define the formal dual (adjoint) problem to
  (\ref{eq:abstract-prob}) to be: For any $g\in L$, find $z(g)\in Z
  \subset V$ such that
  \begin{equation}
    \bid{z(g)}{v} = f_L(g, v) \Foreach v\in V,
  \end{equation}
  where $\bid{v}{w} = \bi{w}{v}$. 

  In addition we suppose that the bilinear form $f_L(\cdot, \cdot)$
  induces a seminorm on $L$, \ie 
  \begin{equation}
    \Norm{\cdot}_L 
    :=
    \sup_{g\in L} \frac{f_L(g, \cdot)}{\Norm{g}_V}.
  \end{equation}
\end{Defn}

\begin{Hyp}
  \label{ass:dual-inf-sup}
  In addition to Assumption \ref{ass:cont-inf-sup}, assume the formal
  dual problem is well posed, that is, for all $g\in L$
  \begin{equation}
    \gamma^*\Norm{z(g)}_Z \leq \sup_{v\in V} \frac{\bid{z(g)}{v}}{\Norm{v}_V},
  \end{equation}
  then we immediately have that
  \begin{equation}
    \Norm{z(g)}_Z \leq \frac{1}{\gamma^*}\Norm{g}_L.
  \end{equation} 
\end{Hyp}

\begin{The}[Dual aposteriori bound]
  \label{the:the-dual-apost}
  Using the same notation as in Theorem
  \ref{the:the-aposteriori-strang} under Assumptions
  \ref{ass:cont-inf-sup} and \ref{ass:dual-inf-sup}
  \begin{equation}
    \label{eq:dual-apost-ab}
    \begin{split}
      \Norm{e}_{L}
      &\leq
      \frac{1}{\gamma^*}
      \bigg(
      \sup_{g\in L, \Norm{g}_L \leq 1}
      \qb{
        l_h(z(g)-z_h) - \bih{u_h}{z(g) - z_h}
      }
      +
      \sup_{g\in L, \Norm{g}_L \leq 1}
      \qb{
        l(z(g)) - l_h(z(g))
      }
      \\ 
      &\qquad \qquad +
      \sup_{g\in L, \Norm{g}_L \leq 1}
      \qb{
        \bih{E(u_h)}{z(g)} - \bi{E(u_h)}{z(g)}
      }
      +
      \sup_{g\in L, \Norm{g}_L \leq 1}
        \bih{E(u_h) - u_h}{z(g)}
      \bigg)
      \\
      & =: 
      \frac{1}{\gamma^*} \qp{ \E_1^L + \E_2^L + \E_3^L + \E_4^L}.
    \end{split}
  \end{equation}
\end{The}
\begin{Proof}
  Proceeding along the same lines as Theorem
  \ref{the:the-aposteriori-strang} we have that $e:= u - u_h = \qp{u -
    E(u_h)} + \qp{E(u_h) - u_h} =: e^C + e^N$ and using the stability
  of the dual problem
  \begin{equation}
    \Norm{e^C}_L
    =
    \sup_{g\in L, \norm{g}_L \leq 1} f_L(e^C, g)
    =
    \sup_{g\in L, \norm{g}_L \leq 1} \bid{z(g)}{e^C}
    =
    \sup_{g\in L, \norm{g}_L \leq 1} \bi{e^C}{z(g)}.
  \end{equation}
  Adding and subtracting appropriate terms gives
  \begin{equation}
    \begin{split}
      \Norm{e^C}_L 
      &=
      \sup_{g\in L, \norm{g}_L \leq 1} 
      l(z(g))
      -
      l_h(z(g))
      +
      l_h(z(g) - z_h)
      -
      \bih{u_h}{z(g) - z_h} 
      \\
      &\qquad\qquad
      -
      \bih{E(u_h) - u_h}{z(g)} + \bih{E(u_h)}{z(g)} -
      \bi{E(u_h)}{z(g)}
      \\
      &\leq
      \sup_{g\in L, \norm{g}_L \leq 1} 
      \bigg[
        l_h(z(g) - z_h)
        -
        \bih{u_h}{z(g) - z_h} 
      \bigg]
+
      \sup_{g\in L, \norm{g}_L \leq 1} 
      \bigg[
        l(z(g))
        -
        l_h(z(g))
      \bigg]
      \\
      &\qquad +
      \sup_{g\in L, \norm{g}_L \leq 1} 
      \bigg[
        \bih{E(u_h)}{z(g)}
        -
        \bi{E(u_h)}{z(g)}
      \bigg]
      +
      \sup_{g\in L, \norm{g}_L \leq 1} 
      \bih{E(u_h) - u_h}{z(g)},
    \end{split}
  \end{equation}
  concluding the proof.
\end{Proof}
\begin{Rem}[Computability of (\ref{eq:dual-apost-ab})]
  As in the primal case Theorem \ref{the:the-dual-apost} is not useful
  in its own right without further assumptions, for example one would
  desire that the finite element space satisfies some approximability
  assumption for any $v\in Z$
  \begin{equation}
    \inf_{v_h \in V_h} \Norm{v - v_h}_{V_h} \leq C_1 h^\alpha \Norm{v}_Z.
  \end{equation}
  The full range of assumptions required for dual a posteriori error
  control tends to be very problem specific as will be illustrated in
  the sequel.
\end{Rem}

\section{Applications to discontinuous Galerkin approximations of 2nd order linear self adjoint problems}
\label{sec:selfadj}

In this section we give an illustrative example showing how to apply
the framework to derive a posteriori bounds for various schemes
approximating a model problem. To that end, suppose that the function
spaces defined in the previous chapter are defined over a domain
$\W$. Let $\T{}$ be a conforming, shape regular triangulation of $\W$,
namely, $\T{}$ is a finite family of sets such that
\begin{enumerate}
\item $K\in\T{}$ implies $K$ is an open simplex (segment for $d=1$,
  triangle for $d=2$, tetrahedron for $d=3$),
\item for any $K,J\in\T{}$ we have that $\closure K\meet\closure J$ is
  a full subsimplex (i.e., it is either $\emptyset$, a vertex, an
  edge, a face, or the whole of $\closure K$ and $\closure J$) of both
  $\closure K$ and $\closure J$ and
\item $\union{K\in\T{}}\closure K=\closure\W$.
\end{enumerate}
The shape regularity of $\T{}$ is defined as the number
\begin{equation}
  \label{eqn:def:shape-regularity}
  \mu(\T{}) := \inf_{K\in\T{}} \frac{\rho_K}{h_K},
\end{equation}
where $\rho_K$ is the radius of the largest ball contained inside
$K$ and $h_K$ is the diameter of $K$. An indexed family of
triangulations $\setof{\T n}_n$ is called \emph{shape regular} if 
\begin{equation}
  \label{eqn:def:family-shape-regularity}
  \mu:=\inf_n\mu(\T n)>0.
\end{equation}
We use the convention where $\funk h\W\reals$ denotes the {piecewise
  constant} \emph{meshsize function} of $\T{}$, i.e.,
\begin{equation}
  h(\vec{x}):=\max_{\closure K\ni \vec x}h_K,
\end{equation}
which we shall commonly refer to as $h$.

We let $\E{}$ be the skeleton (set of common interfaces) of the
triangulation $\T{}$ and say $e\in\E$ if $e$ is on the interior of
$\W$ and $e\in\partial\W$ if $e$ lies on the boundary $\partial\W$
{and set $h_e$ to be the diameter of $e$.}


We let $\poly k(\T{})$ denote the space of piecewise polynomials of
degree $k$ over the triangulation $\T{}$,\ie
\begin{equation}
  \poly k (\T{}) = \{ \phi \text{ such that } \phi|_K \in \poly k (K) \}
\end{equation}
 and introduce the \emph{finite element space}
\begin{gather}
  \label{eqn:def:finite-element-space}
  \fes := \dg\qp{k} = \poly k(\T{}) 
\end{gather}
to be the usual space of discontinuous piecewise polynomial
functions of degree $k$.

\begin{Defn}[broken Sobolev spaces, trace spaces]
  \label{defn:broken-sobolev-space}
  We introduce the broken Sobolev space
  \begin{equation}
    \sobh{l}(\T{})
    :=
    \ensemble{\phi}
             {\phi|_K\in\sobh{l}(K), \text{ for each } K \in \T{}}.
  \end{equation}
  We also make use of functions defined in these broken spaces
  restricted to the skeleton of the triangulation. This requires an
  appropriate trace space
  \begin{equation}
    \Tr{\E} := \prod_{K\in\T{}} \leb{2}(\partial K).
  \end{equation}
\end{Defn}

\begin{Defn}[jumps, averages and tensor jumps]
  \label{defn:averages-and-jumps}
  We may define average, jump and tensor jump operators over $\Tr{\E}$ for
  arbitrary scalar functions $v\in\Tr{\E}$ and vectors $\vec v\in\Tr{\E}^d$.
  \begin{gather*}
    \label{eqn:average}
    \dfunkmapsto[]
	        {\avg{\cdot}}
	        v
	        {\Tr{\E}}
	        {
                  \begin{cases}
                    \frac{1}{2}\qp{v|_{K_1} + v|_{K_2}} 
                    \\
                    v|_{\partial\W} \text{ on } \partial\W
                  \end{cases}
                }
	        {\leb{2}(\E)}
                \qquad\qquad
                \dfunkmapsto[]
	        {\avg{\cdot}}
	        {\vec v}
	        {\qb{\Tr{\E}}^d}
	        {
                  \begin{cases}
                    \frac{1}{2}\qp{\vec{v}|_{K_1} + \vec{v}|_{K_2}}
                    \\
                    \vec v|_{\partial\W} \text{ on } \partial\W
                  \end{cases}
                }
	        {\qb{\leb{2}(\E\cup \partial\W)}^d}
                \\
                \dfunkmapsto[]
	        {\jump{\cdot}}
	        {v}
	        {{\Tr{\E}}}
	        {
                  \begin{cases}
                    {{v}|_{K_1} \geovec n_{K_1} + {v}|_{K_2}} \geovec n_{K_2}
                    \\
                    \qp{v \vec n}|_{\partial\W} \text{ on } \partial\W
                  \end{cases}
                }
	        {\qb{\leb{2}(\E)}^d}
                                \qquad\qquad
    \dfunkmapsto[]
	        {\jump{\cdot}}
	        {\vec v}
	        {\qb{\Tr{\E}}^d}
	        {
                  \begin{cases}
                    {\Transpose{\qp{\vec{v}|_{K_1}}}\geovec n_{K_1} 
                      +
                      \Transpose{\qp{\vec{v}|_{K_2}}}\geovec n_{K_2}}
                    \\
                    \qp{\Transpose{\vec v} \vec n}|_{\partial\W} \text{ on } \partial\W
                  \end{cases}
                }
	        {{\leb{2}(\E)}}
\\
                \label{eqn:tensor-jump}
    \dfunkmapsto[.]
	        {\tjump{\cdot}}
	        {\vec v}
	        {\qb{\Tr{\E\cup \partial\W}}^d}
	        {
                  \begin{cases}
                    {\vec{v}|_{K_1} }\otimes \geovec n_{K_1} + \vec{v}|_{K_2} \otimes\geovec n_{K_2}
                    \\
                    \qp{{\vec v} \otimes \vec n}|_{\partial\W} \text{ on } \partial\W
                  \end{cases}
                }
	        {\qb{\leb{2}(\E\cup \partial\W)}^{d\times d}}
  \end{gather*}
\end{Defn}

\subsection{Conforming reconstruction
  operators}

A simple, quite general methodology for the construction of
reconstruction operators with the desirable properties mentioned in
Remark \ref{rem:computability} is to use an averaging interpolation
operator into an $\sobh{s}$ ($s=1,2$) conforming finite element
space. For example the Oswald interpolant into a $\cont{0}$ Lagrangian
finite element space can be used for a $\sobh{1}$ conforming operator
\cite{KarakashianPascal:2003} or a $\cont{1}$ Hsieh--Clough--Tocher
macro element conforming space for $\sobh{2}$ conformity
\cite[c.f.]{BrennerGudiSung:2010,GeorgoulisHoustonVirtanen:2011}.

\begin{Example}[Lowest order $\sobh{1}(\W)$ and $\sobh{2}(\W)$ reconstructions]
  An example of the lowest order ($p=1$) $\sobh{1}(\W)$ reconstruction
  operator $E^1(u_h)$ is for given $u_h\in\fes$ to specify the values
  of the reconstruction at the vertices as an average in a local
  neighbourhood, as with the Oswald interpolant \cite{ErnGuermond:2004}. In general let $\vec x$ be a degree of freedom of the
  conformal space $\fesW$ and let $\patch{K_{\vec x}}$ be the
  set of all elements sharing the degree of freedom $\vec x$ then the
  reconstruction at that specific degree of freedom is given by
  \begin{equation}
    E(u_h)(\vec x) 
    =
    \frac{1}{\text{card}(\patch{K_{\vec x}})}\sum_{K \in \patch{K_{\vec x}}}  u_h\vert_K(\vec x).
  \end{equation}
  In the case $k=1$ the associated degrees of freedom are illustrated below:
  \begin{center}
    \begin{tikzpicture}      
     
      \path[fill=green!60] 
      (0,0) node (N0) {}
      -- (2,0) node (N5) {}
      -- (2,2) node (N2) {}
      -- cycle;
      
      \path[coordinate]
      (2,.67) node (K) {};
      
      \path[fill=green!60,xshift=6cm] 
      (0,0) node (NL1) {}
      -- (2,0) node (NL2) {}
      -- (2,2) node (NL0) {}
      -- cycle;
      
      \foreach \s in {0,1}{
	\path[coordinate,xshift=6cm*\s]
	(1.33,.67) node (L\s) {};
      }
      
      \path[fill=green!60]
      (4,0) node (N1) {}
      -- (2,2)
      -- (2,0)
      -- cycle
      ;
      
      \path[fill=green!60,xshift=6cm]
      (4,0) node (NR0) {}
      -- (2,2) node (NR1) {}
      -- (2,0) node (NR2) {}
      -- cycle
      ;
      
      \foreach \s in {0,1}{
	\path[coordinate,xshift=6cm*\s]
	(2.67,.67) node (R\s) {};
      }
      
      \path[coordinate]
      (3,1) node (N3) {}
      (1,1) node (N4) {}
      ;
      
      \path[coordinate,xshift=6cm]
      (1,0) node (NL3) {}
      (2,1) node (NL4) {}
      (1,1) node (NL5) {}
      ;

      \path[coordinate,xshift=6cm]
      (2,1) node (NR3) {}
      (3,0) node (NR4) {}
      (3,1) node (NR5) {}
      ;      
      
      \foreach \s in {0,1,...,2}{
	\node (K\s) at (N\s) [circle,inner sep=1pt,fill=blue!30,draw] {\phantom{1}};
      }
      
      \foreach \s in {0,1}{
	\node (KL\s) at (NL\s) [circle,inner sep=1pt,fill=blue!30,draw] {\phantom{1}};
	\node (KR\s) at (NR\s) [circle,inner sep=1pt,fill=blue!30,draw] {\phantom{1}};
      }
      
      
      \node at (K) [rectangle,inner sep=1pt] {$u_h$};
      \node at (L1) [rectangle,inner sep=1pt,xshift=.6cm ] {$E^1(u_h)$};

      \foreach \s in {0,1}{
      }
      
      
      \draw[->] (4.5,0.9) -- (5.5,.9)
      node[pos=0.5,above] {$\sobh1(\W)$ reconstruction};
    \end{tikzpicture}
  \end{center}
  
  Notice that the degrees of freedom of the reconstruction match that
  of the original function, as such $E^1(u_h)\in\fes$, this is because
  of the existence of a conforming $\sobh1(\W)$ subspace in
  $\fes$. For $\sobh{1}(\W)$ conformity this is true for any $k$, for
  $\sobh2(\W)$ one would need $k \geq 5$ to have access to appropriate
  elements, as such a richer space is required. For $k=2$ the
  $\poly{4}$ Hsieh--Clough--Tocher macro element is used with degrees
  of freedom given below:
  \begin{center}
    \begin{tikzpicture}      
      
      \path[fill=green!60] 
      (0,0) node (N0) {}
      -- (2,0) node (N5) {}
      -- (2,2) node (N2) {}
      -- cycle;
      
      \path[coordinate]
      (2,.67) node (K) {};
      
      \path[fill=green!60,xshift=6cm] 
      (0,0) node (NL1) {}
      -- (2,0) node (NL2) {}
      -- (2,2) node (NL0) {}
      -- cycle;
      
      \foreach \s in {0,1}{
	\path[coordinate,xshift=6cm*\s]
	(1.33,.67) node (L\s) {};
      }
      
      \path[fill=green!60]
      (4,0) node (N1) {}
      -- (2,2)
      -- (2,0)
      -- cycle
      ;
      
      \path[fill=green!60,xshift=6cm]
      (4,0) node (NR0) {}
      -- (2,2) node (NR1) {}
      -- (2,0) node (NR2) {}
      -- cycle
      ;
      
      \foreach \s in {0,1}{
	\path[coordinate,xshift=6cm*\s]
	(2.67,.67) node (R\s) {};
      }
      
      \path[coordinate]
      (3,1) node (N3) {}
      (1,1) node (N4) {}
      ;
      
      \path[coordinate,xshift=6cm]
      (1,0) node (NL3) {}
      (2,1) node (NL4) {}
      (1,1) node (NL5) {}
      ;

      \path[coordinate,xshift=6cm]
      (2,1) node (NR3) {}
      (3,0) node (NR4) {}
      (3,1) node (NR5) {}
      ;      
      
      \foreach \s in {0,1,...,5}{
	\node (K\s) at (N\s) [circle,inner sep=1pt,fill=blue!30,draw] {\phantom{1}};
      }

      \foreach \s in {0,1,2,5}{
	\node (KL\s) at (NL\s) [circle,inner sep=1pt,fill=blue!30,draw] {\phantom{1}};
	\node (KR\s) at (NR\s) [circle,inner sep=1pt,fill=blue!30,draw] {\phantom{1}};
      }
      \foreach \s in {0,1}{
	\node (KL\s) at (NL\s) [circle,inner sep=1pt,fill=red!30,draw] {};
	\node (KR\s) at (NR\s) [circle,inner sep=1pt,fill=red!30,draw] {};
      }

      \path[coordinate,xshift=6cm]
      (2,1) node (C1) {}
      ;      
      
      \node (K1) at (C1) [circle,inner sep=1pt,fill=blue!30,draw] {\phantom{1}};
      \node (K1) at (C1) [circle,inner sep=1pt,fill=red!30,draw] {};

      \draw[->] (7,0) -- (7,-0.5);
      \draw[->] (9,0) -- (9,-0.5);
      
      \draw[->] (6.5,0.5) -- (6,1);
      \draw[->] (7.5,1.5) -- (7,2);

      \draw[->] (9.5,0.5) -- (10,1);
      \draw[->] (8.5,1.5) -- (9,2);

      
      \node at (K) [rectangle,inner sep=1pt] {$u_h$};
      \node at (L1) [rectangle,inner sep=1pt,xshift=.6cm ] {$E^2(u_h)$};

      \foreach \s in {0,1}{
      }
      
      
      \draw[->] (4.5,0.9) -- (5.5,0.9)
      node[pos=0.5,above] {$\sobh2(\W)$ reconstruction};
    \end{tikzpicture}
\end{center}
Note that the reconstruction is no longer an element of $\fes$.
\end{Example}


\begin{Defn}[$\sobh1$ and $\sobh2$ mesh dependant norms]
  We define the $\sobh1$ and $\sobh2$ mesh dependant norms to be 
  \begin{gather}
    \enorm{w_h}{2}^2 := \Norm{\nabla_h w_h}_{\leb{2}(\W)}^2 + h^{-1}_e\Norm{\jump{u_h}}_{\leb{2}(\E)}^2
    \\
    \eenorm{w_h}{2}^2 
    :=
    \Norm{\Hess_h w_h}_{\leb{2}(\W)}^2
    +
    h_e^{-1} \Norm{\jump{\nabla w_h}}_{\leb{2}(\E)}^2
    +
    h_e^{-3} \Norm{\jump{w_h}}_{\leb{2}(\E)}^2.
  \end{gather}
\end{Defn}

\begin{Lem}[Reconstruction bounds {\cite{KarakashianPascal:2003,GeorgoulisHoustonVirtanen:2011}}]
  \label{lem:reconstruction-bounds}
  There exist operators $E^s : \fes\to\sobh{s}(\W)$, $s=1,2$ such that 
  \begin{gather}
      \Norm{E^1(u_h) - u_h}_{\leb{2}(\W)}^2 
      \leq
      C h_e \Norm{\jump{u_h}}_{\leb{2}(\E)}^2
      \\
      \enorm{E^1(u_h) - u_h}{2}^2 
      \leq C h_e^{-1}
      \Norm{\jump{u_h}}_{\leb{2}(\E)}^2
  \end{gather}
  and
  \begin{equation}
    \eenorm{E^2(u_h) - u_h}{2}^2
    \leq
    C \qp{
      h_e^{-1} \Norm{\jump{\nabla u_h}}_{\leb{2}(\E)}^2
      + 
      h_e^{-3} \Norm{\jump{u_h}}_{\leb{2}(\E)}^2
    }.
  \end{equation}
\end{Lem}
\begin{Proof}
  The proof of the bound for $E^1(u_h)$ is given in \cite[Thm 2.2]{KarakashianPascal:2003} and for $E^2(u_h)$ is given in \cite[Lem 3.1]{GeorgoulisHoustonVirtanen:2011}.
\end{Proof}

\subsection{Laplace's problem}

We begin by setting $V = W = \hoz(\W)$, which is the subspace of $\sobh1(\W)$ with functions vanishing on $\partial\W$, and let
\begin{equation}
  \bi{u}{v} = \int_\W \nabla u \cdot \nabla v.
\end{equation}
This operator is coercive and thus certainly inf-sup stable in $\hoz(\W)$, thus satisfies Assumption \ref{ass:cont-inf-sup}. We consider discontinuous Galerkin approximations with $V_h = W_h = \poly{k}(\T{})$. 

The interior penalty (IP) method \cite{Arnold:1982}
\begin{equation}
  \label{eq:IP}
  \bih{u_h}{v_h} 
  =
  \int_\W \nabla_h u_h \cdot \nabla_h v_h
  -
  \int_\E \qp{\jump{u_h} \cdot \avg{\nabla v_h}
  +
  \jump{v_h} \cdot \avg{\nabla u_h}
  -
  \frac{\sigma}{h} \jump{u_h}\cdot \jump{v_h}}
\end{equation}
is well known to be bounded in the $\sobh{1}(\W)$ like dG--norm
$\enorm{\cdot}{2}$. A conforming space $\fesW\subset \hoz(\W)$ can be
found by enforcing continuity of the piecewise polynomial functions,
\ie taking $\fesW = \poly{p}(\T{})\cap \hoz(\W)$. Reconstruction operators
$E$ are not unique in this setting. Indeed, one example can be found
in \cite{KarakashianPascal:2003} and another the Ritz projection of
$u_h$ into $\poly{p}(\T{})\cap \hoz(\W)$. 

\begin{Rem}[Artificial regularity requirement]
  The IP method given in (\ref{eq:IP}) is only consistent for
  $u\in\sobh{3/2+\epsilon}(\W)$ since we require $\avg{\nabla u}$ to
  be well defined in order to invoke Strang's Lemma. This regularity
  is all that is required for the a posteriori analysis based on
  coercivity given in \cite{KarakashianPascal:2003} for example but is
  insufficient for the abstract framework devleoped in
  \S\ref{sec:setup}, for this we would also require
  $v\in\sobh{3/2+\epsilon}(\W)$ (in the evaluation of $\bih{u_h}{v -
    v_h}$ for example) which is incompatible with the stability
  theory. For the bound in Theorem \ref{the:the-aposteriori-strang} to
  make sense it is required that the discrete bilinear form
  $\bih{\cdot}{\cdot}$, which is well defined over $\sobh2{(\T{})}
  \times \sobh2{(\T{})}$, is extended to $\sobh{1}(\T{}) \times
  \sobh{1}(\T{})$, the correct domain to invoke the appropriate
  stability theory used for this class of PDE.
  
  The extension is not unique. Indeed, one method of extending is to
  modify the definition of the bilinear form for $u,v
  \in\sobh{1}(\T{})$ to
  \begin{equation}
    \label{eq:IP-mod}
    \bih{u}{v} 
    =
    \int_\W \nabla_h u \cdot \nabla_h v
    -
    \int_\E \qp{\jump{u} \cdot \avg{P_{k-1} \qp{\nabla v}}
      +
      \jump{v} \cdot \avg{P_{k-1} \qp{\nabla u}}
      -
      \frac{\sigma}{h} \jump{u}\cdot \jump{v}},
  \end{equation}
  where $P_{k-1}$ denotes the $\leb{2}$ projection into
  $\poly{k-1}(\T{})$. Note that this is very similar to a proceedure
  proposed in \cite[\S 4]{AAFCS:2009}, however extending the bilinear
  form in this fashion means that (\ref{eq:IP-mod}) coincides with
  (\ref{eq:IP}) over $\fes\times \fes$, as such (\ref{eq:IP-mod})
  \emph{is} the classical IP method interpreted in such a way that a
  posteriori analysis can be conducted using the tools of
  \S\ref{sec:setup}. 

  By modifying the bilinear form in this fashion the method becomes
  inconsistent (in the a priori sense of (\ref{eq:apriori-strang})),
  but the inconsistency term is controllable and convergence is
  optimal. Another method to see this is given in
  \cite{Gudi:2010}. Here the author uses a posteriori techniques and
  the averaging operator from \cite{KarakashianPascal:2003} to control
  the inconsistency arising from the lack of regularity up to data
  oscillation terms.  The a posteriori bound derived under the
  framework of \S\ref{sec:setup} is not effected by the inconsistency
  as we will see in the sequel.
\end{Rem}

As the IP method, the Babu{\v{s}}hka--Zl\'amal (BZ) method
\cite{BabuskaZlamal:1973}
\begin{equation}
  \bih{u_h}{v_h} 
  =
  \int_\W \nabla_h u_h \cdot \nabla_h v_h
  +
  \int_\E
  \frac{\sigma}{h} \jump{u_h}\cdot \jump{v_h}
\end{equation}
is bounded in the dG--norm, $\enorm{\cdot}{2}$, however it is
inconsistent regardless of the regularity of $u$, that is the second
term in (\ref{eq:apriori-strang}) is always nonzero. It is of optimal
order hence, in view of Strang's Lemma, the method converges
optimally. Note also that this bilinear form can be trivially extended
to $\sobh{1}(\T{}) \times \sobh{1}(\T{})$.

\begin{Pro}[Approximation properties of $\leb{2}$ projections {\cite[{c.f. \S 5.6.2.2}]{Di-PietroErn:2012}}]
  \label{pro:approx-l2proj}
  Let $P_0$ denote the $\leb{2}$ orthogonal projection into
  $\poly{0}(\T{})$ then for any $v\in\sobh{1}(\T{})$ we have the
  following local approximation bounds over elements, $K\in\T{}$ and
  faces $e\in\E$:
  \begin{gather}
    \Norm{v- P_0 v}_{\leb{2}(K)} \leq C h_K \Norm{\nabla v}_{\leb{2}(K)}
    \\
    \Norm{v- P_0 v}_{\leb{2}(e)} \leq C h_e^{1/2} \Norm{\nabla v}_{\leb{2}(K)}.
  \end{gather}
\end{Pro}

\begin{Pro}[Approximation properties of the Scott-Zhang interpolator {\cite[{c.f. \S 1.130}]{ErnGuermond:2004}}]
  \label{pro:SZ}
  Let $I_k$ denote the Scott-Zhang interpolant into
  $\fes\cap\hoz(\W)$, then the following local approximation bounds
  over elements, $K\in\T{}$ and faces $e\in\E$ hold:
  \begin{gather}
    \Norm{v- I_k v}_{\leb{2}(K)} \leq C h_K \Norm{\nabla v}_{\leb{2}(\patch{K})}
    \\
    \Norm{v- I_k v}_{\leb{2}(e)} \leq C h_e^{1/2} \Norm{\nabla v}_{\leb{2}(\patch{K})},
  \end{gather}
  where $\widehat{K}$ denotes a \emph{patch} of an element $K$.
\end{Pro}

\begin{Lem}[Energy norm a posteriori control]
  \label{lem:lap-apost}
  Let $u_h$ be either the IP or BZ approximation of the solution of
  Laplace's problem, then up to a constant independent of the meshsize
  both approximations are a posteriori controllable by the following bound:
      \begin{equation}
        \begin{split}
          \enorm{u - u_h}{2} 
          &\leq
          C 
          \qp{\sum_{K\in\T{}} \eta_{1,R}^2 + \sum_{e \in \E} \eta_{1,J}^2}^{1/2} \AND
        \end{split}        
      \end{equation}
      with 
      \begin{equation}
        \begin{split}
          \eta_{1,R}^2 
          &:=
          h_K^2 \Norm{f + \Delta u_h}_{\leb{2}(K)}^2
          \\
          \eta_{1,J}^2
          &:=
          h_e \Norm{\jump{\nabla u_h}}_{\leb{2}(e)}^2
          +
          \sigma h_e^{-1} \Norm{\jump{u_h}}_{\leb{2}(e)}^2.
        \end{split}
      \end{equation}
\end{Lem}
\begin{Proof}
  Apply Theorem \ref{the:the-aposteriori-strang}, noting the
  inconsistency terms for \emph{both} methods vanish. Then use
  standard a posteriori arguments
  \cite[c.f.]{Ainsworth:2000,KarakashianPascal:2003} for the residual
  term which we summarise here for completeness. Firstly, for the IP
  method
  \begin{equation}
    \label{eq:pf1}
    \begin{split}
    l_h(v-v_h) - \bih{u_h}{v - v_h}
    &=
    \int_\W f\qp{v - v_h} - \nabla u_h \cdot \qp{\nabla v - \nabla v_h}
    +
    \int_\E \jump{u_h} \cdot \avg{P_{k-1}\qp{\nabla v - \nabla v_h}}
    \\&\qquad 
    + 
    \int_\E
    \jump{v -v_h} \avg{\nabla u_h}
    -
    \frac{\sigma}{h} \jump{u_h} \cdot \jump{{v - v_h}}
    \\
    &=
    \int_\W \qp{f + \Delta u_h} \qp{v - v_h} - \int_\E \jump{\nabla u_h} \avg{v - v_h}
    \\
    & \qquad + 
    \int_\E \jump{u_h} \cdot \avg{P_{k-1}\qp{\nabla v - \nabla v_h}}
    -
    \frac{\sigma}{h} \jump{u_h} \cdot \jump{{v - v_h}}
    \\
    &=: \cR_1 + \cR_2 + \cR_3 + \cR_4.
  \end{split}
  \end{equation}
  Using Cauchy-Schwarz we see
  \begin{equation}
    \label{eq:pf2}
    \begin{split}
      \cR_1 
      &=
      \int_\W \qp{f + \Delta u_h} \qp{v - v_h} 
      \\
      &=
      \sum_{K\in\T{}}
      \Norm{f + \Delta u_h}_{\leb{2}(K)} 
      \Norm{v - v_h}_{\leb{2}(K)} 
      \\
      &\leq
      C
      \Norm{\nabla v}_{\leb{2}(\W)} 
      \sum_{K\in\T{}}
      h_K
      \Norm{f + \Delta u_h}_{\leb{2}(K)},
    \end{split}
  \end{equation}
  by taking $v_h = P_0 v$, the $\leb{2}$ orthogonal projection into
  $\poly{0}(\T{})$ and using Proposition
  \ref{pro:approx-l2proj}. For the second term, using the same $v_h$
  \begin{equation}
    \label{eq:pf3}
    \begin{split}
    \cR_2
    &=
    -\sum_{e\in\E} \int_e \jump{\nabla u_h} \avg{v - v_h}
    \\
    &\leq
    \sum_{e\in\E} \Norm{\jump{\nabla u_h}}_{\leb{2}(e)} \Norm{ \avg{v - v_h}}_{\leb{2}(e)}
    \\
    &\leq
    C \Norm{\nabla v}_{\leb{2}(\W)}
    \sum_{e\in\E} h_e^{1/2} \Norm{\jump{\nabla u_h}}_{\leb{2}(e)}.
    \end{split}
  \end{equation}
  For the third term, in view of Cauchy Schwarz and a trace inequality
  and the stability of the $\leb{2}(\W)$ projection
  \begin{equation}
    \label{eq:pf4}
    \begin{split}
      \cR_3
      &=
      \sum_{e\in\E}
      \int_e \jump{u_h} \cdot \avg{{P_{k-1}\qp{\nabla v - \nabla v_h}}}
      \\
      &\leq
      C
      \sum_{e\in\E}
      \Norm{\jump{u_h}}_{\leb{2}(e)} \Norm{ \avg{P_{k-1}\qp{\nabla v - \nabla v_h}}}_{\leb{2}(e)}
      \\
      &\leq
      C
      \Norm{ \nabla v }_{\leb{2}(\W)}
      \sum_{e\in\E}
      h_e^{-1/2}
      \Norm{\jump{u_h}}_{\leb{2}(e)},
    \end{split}
  \end{equation}
  and similarly for the final term
  \begin{equation}
    \label{eq:pf5}
    \begin{split}
      \cR_4 
      &=
      \sum_{e\in\E}
      \int_e \sigma h_e^{-1} \jump{u_h} \cdot \jump{{v - v_h}}
      \\
      &\leq 
      \sum_{e\in\E}
      \int_e
      \sigma h_e^{-1} \Norm{\jump{u_h}}_{\leb{2}(e)} \Norm{\jump{{v - v_h}}}_{\leb{2}(e)}
      \\
      &\leq
      C\sigma
      \Norm{ \nabla v }_{\leb{2}(\W)}
      \sum_{e\in\E}
      h_e^{-1/2} \Norm{\jump{u_h}}_{\leb{2}(e)}.
    \end{split}
  \end{equation}
  Collecting (\ref{eq:pf1})--(\ref{eq:pf5}) we have
  \begin{equation}
    \begin{split}
      \cE_1 
      &:= 
      \sup_{v\in \hoz, \Norm{\nabla v}_{\leb{2}(\W)} \leq 1}
      \qb{
        l_h(v-v_h) - \bih{u_h}{v - v_h}
      }
      \\
      &\leq
      C \qp{\sum_{K\in\T{}} \eta_{1,R}^2 + \sum_{e \in \E} \eta_{1,J}^2}^{1/2} 
    \end{split}
  \end{equation}
  yielding the desired result. For the BZ method we note
  \begin{equation}
    \label{eq:pf1}
    \begin{split}
      l_h(v-v_h) - \bih{u_h}{v - v_h}
      &=
      \int_\W f\qp{v - v_h} - \nabla u_h \cdot \qp{\nabla v - \nabla v_h}
      -
      \frac{\sigma}{h} \jump{u_h} \cdot \jump{{\nabla v - \nabla v_h}}
      \\
      &=
      \int_\W \qp{f + \Delta u_h} \qp{v - v_h} - \int_\E \jump{\nabla u_h} \avg{v - v_h}
      \\
      & \qquad +
      \int_\E \jump{v - v_h} \cdot \avg{{\nabla u_h}}
      -
      \frac{\sigma}{h} \jump{u_h} \cdot \jump{{v - v_h}}
      \\
      &=: \cR_1 + \cR_2 + \widetilde{\cR_3} + \cR_4.
    \end{split}
  \end{equation}
  Notice that $\cR_1, \cR_2 \AND \cR_4$ are the same terms that appear
  in the analysis of the IP method above. The third term is different
  due to the nature of the inconsistency of the scheme. This term vanishes as
  soon as $v_h$ is chosen as a continuous function. It suffices
  therefore to pick $v_h$ as the Scott--Zhang interpolator into
  $\fes\cap \hoz(\W)$ and making use of Proposition \ref{pro:SZ}. The
  bounds for the other terms follow analogously as above, concluding
  the proof.
\end{Proof}

\begin{Rem}[Lower bounds]
  It can be proven that the a posteriori bounds for both IP and BZ
  methods satifsy lower bounds to the error, that is, using the
  notation of Lemma \ref{lem:lap-apost}
  \begin{equation}
    \label{eq:lower-bounds}
    \qp{\sum_{K\in\T{}} \eta_{1,R}^2 + \sum_{e \in \E} \eta_{1,J}^2}^{1/2}
    \leq
    \enorm{u - u_h}{2} 
    + 
    h_K\Norm{f - P_k f}_{\leb{2}(\W)}.
  \end{equation}
  The last term in (\ref{eq:lower-bounds}) is called a \emph{data
    oscillation} term. See \cite{KarakashianPascal:2007} for a proof
  of this fact for the IP method. The arguement for the BZ method is
  identical since the a posteriori bound is the same. We will examine
  how other forms of inconsistency effect lower bounds in the sequel.
\end{Rem}

\begin{Rem}[Dual norm estimates]
  Obtaining optimal dual norm estimates is slightly more delicate than
  the primal ones. Suppose the method is consistent, and in particular
  adjoint consistent \cite{ArnoldBrezziCockburnMarini:2001}, that is
  if $z\in\sobh{2}(\W)$ solves the dual problem $-\Delta z = g$ then
  \begin{equation*}
    \label{eq:adjointconsistent}
    \bih{v}{z} = \int_\W v g \Foreach v\in\sobh2\qp{\T{}}.
  \end{equation*}
  The IP method satisfies this condition and optimal a posteriori bounds can be obtained in the dual norm quite simply under the framework of Theorem \ref{the:the-dual-apost} as the following Lemma illustrates:
\end{Rem}

\begin{Lem}[Adjoint consistent dual norm a posteriori control]
  \label{lem:lap-apost-dual}
  Let $u_h$ be the IP approximation of the solution of
  Laplace's problem, then up to a constant independent of the meshsize
  the following a posteriori bound holds:
      \begin{equation}
        \begin{split}
          \Norm{u - u_h}_{\leb{2}(\W)}
          &\leq
          C 
          \qp{\sum_{K\in\T{}} \eta_{0,R}^2 + \sum_{e \in \E} \eta_{0,J}^2}^{1/2} 
        \end{split}        
      \end{equation}
      with 
      \begin{equation}
        \begin{split}
          \eta_{0,R}^2 
          &:=
          h_K^4 \Norm{f + \Delta u_h}_{\leb{2}(K)}^2
          \\
          \eta_{0,J}^2
          &:=
          h_e^3 \Norm{\jump{\nabla u_h}}_{\leb{2}(e)}^2
          +
          \sigma h_e \Norm{\jump{u_h}}_{\leb{2}(e)}^2.
        \end{split}
      \end{equation}
\end{Lem}
\begin{Proof}
  To control the nonconformity term
  arising in the dual a posteriori bound ($\cE_4^L$ in
  (\ref{eq:dual-apost-ab})) we note
  \begin{equation}
    \begin{split}
      \bih{E^1(u_h) - u_h}{z} 
      &=
      \int_\W \nabla_h\qp{E^1(u_h) - u_h} \cdot \nabla z
      +
      \int_\E \jump{u_h} \cdot \avg{\nabla z}
      \\
      &=
      \int_\W -\qp{E^1(u_h) - u_h} \Delta z.
    \end{split}
  \end{equation}
  Hence 
  \begin{equation}
    \begin{split}
      \cE_5 
      &=
      \sup_{g\in L, \norm{g}_L \leq 1} \bih{E^1(u_h) - u_h}{z}
      \\
      &=
      \sup_{g\in \leb{2}(\W), \Norm{g}_{\leb{2}(\W)} \leq 1}
      \int_\W -\qp{E^1(u_h) - u_h} \Delta z
      \\
      &\leq
      \sup_{g\in \leb{2}(\W), \Norm{g}_{\leb{2}(\W)} \leq 1}
      \Norm{E^1(u_h) - u_h}_{\leb{2}(\W)} \Norm{\Delta z}_{\leb{2}(\W)}
      \\
      &\leq
      \frac{1}{\gamma*} \Norm{E^1(u_h) - u_h}_{\leb{2}(\W)}
      \\
      &\leq
      C h_e^{1/2} \Norm{\jump{u_h}}_{\leb{2}(\E)}.
    \end{split}
  \end{equation}
  The bounds for the other terms follow from standard a posteriori
  arguments and can be derived by mimicking the proof of Lemma
  \ref{lem:lap-apost} and using regularity bounds for the solution of
  the dual problem.
\end{Proof}

\begin{Rem}[Dual bounds do not require reconstructions]
  To construct the dual a posteriori bound in $\leb{2}$ for adjoint
  consistent methods it is not necessary to use the reconstruction
  approach. The numerical approximation is smooth enough to directly
  apply the stability theory of the dual problem directly, at least in
  the case of Laplacian and biharmonic operators
  \cite{GeorgoulisVirtanen:2013}.
\end{Rem}

\begin{Rem}[Dual norm estimates for inconsistent methods]
  For adjoint inconsistent methods, i.e., those not satisfying the
  condition in Remark \ref{eq:adjointconsistent}, to obtain optimal
  dual norm a posteriori estimates over-penalisation is necessary, as
  in the a priori case \cite[\S 5.2]{ArnoldBrezziCockburnMarini:2001}. We consider the over-penalised BZ method with
  \begin{equation}
    \label{eq:BZ-overpen}
    \bih{u_h}{v_h}
    =
    \int_\W
    \nabla_h u_h \cdot \nabla_h v_h
    +
    \int_\E
    \sigma h_e^{-\beta} \jump{u_h} \cdot \jump{v_h},
  \end{equation}
  for $\beta \geq 1$ and define the over-penalised norm as 
  \begin{equation}
    \label{eq:opdgnorm}
    \opnorm{u_h}{2} := \Norm{\nabla_h w_h}_{\leb{2}(\W)}^2 + h^{-\beta}_e\Norm{\jump{u_h}}_{\leb{2}(\E)}^2.
  \end{equation}
  Optimal a posteriori bounds can be obtained for the over-penalised BZ method with $\beta$ large enough under the framework of Theorem \ref{the:the-dual-apost} as the following Lemma illustrates:
\end{Rem}

\begin{Lem}[Adjoint inconsistent dual norm a posteriori control]
  \label{lem:lap-apost-dual-incons}
  Let $u_h$ be the over-penalised BZ approximation of the solution of
  Laplace's problem, then for $\beta \geq 3$ up to a constant independent of the meshsize the following a posteriori bound holds:
      \begin{equation}
        \begin{split}
          \Norm{u - u_h}_{\leb{2}(\W)}
          &\leq
          C 
          \qp{\sum_{K\in\T{}} \eta_{0,R}^2 + \sum_{e \in \E} \eta_{0,J}^2}^{1/2} 
        \end{split}        
      \end{equation}
      with 
      \begin{equation}
        \begin{split}
          \eta_{0,R}^2 
          &:=
          h_K^4 \Norm{f + \Delta u_h}_{\leb{2}(K)}^2
          \\
          \eta_{0,J}^2
          &:=
          h_e^3 \Norm{\jump{\nabla u_h}}_{\leb{2}(e)}^2
          +
          \sigma h_e \Norm{\jump{u_h}}_{\leb{2}(e)}^2.
        \end{split}
      \end{equation}
\end{Lem}
\begin{Proof}
  The proof is analogous to Lemma \ref{lem:lap-apost-dual} with the exception of the nonconforming term, here
  \begin{equation}
    \begin{split}
      \bih{E^1(u_h) - u_h}{z} 
      &=
      \int_\W \nabla_h\qp{E^1(u_h) - u_h} \cdot \nabla z
      \\
      &=
      \int_\W -\qp{E^1(u_h) - u_h} \Delta z.
      +
      \int_\E \jump{E^1(u_h) - u_h} \cdot \avg{\nabla z}
    \end{split}
  \end{equation}
  Following the arguments of \cite[\S
  5.2]{ArnoldBrezziCockburnMarini:2001} we have that
  \begin{equation}
    \begin{split}
      \int_\E \jump{E^1(u_h) -u_h} \cdot \avg{\nabla z}
      &\leq
      C h_e^{-\beta/2} \Norm{\jump{E^1(u_h) -u_h}}_{\leb{2}(\E)}  h_e^{\beta/2} \Norm{\avg{\nabla z}}_{\leb{2}(\E)} 
      \\
      &\leq
      C \opnorm{E^1(u_h) -u_h}{2} h_e^{\qp{\beta-1}/{2}} \Norm{\Delta z}_{\leb{2}(\W)}
    \end{split}
  \end{equation}
  via a trace inequality and the definition of the over-penalised dG
  norm (\ref{eq:opdgnorm}). Now through an inverse inequality we see
  \begin{equation}
    \begin{split}
      \int_\E \jump{E^1(u_h) -u_h} \cdot \avg{\nabla z}
      \leq
      C h_e^{\qp{\beta-3}/{2}} \Norm{E^1(u_h) -u_h}_{\leb{2}(\W)}  \Norm{\Delta z}_{\leb{2}(\W)},
    \end{split}
  \end{equation}
  and hence
  \begin{equation}
    \begin{split}
      \cE_5 
      &=
      \sup_{g\in L, \norm{g}_L \leq 1} \bih{E^1(u_h) - u_h}{z}
      \\
      &=
      \sup_{g\in \leb{2}(\W), \Norm{g}_{\leb{2}(\W)} \leq 1}
      \qb{\int_\W -\qp{E^1(u_h) - u_h} \Delta z  +
      \int_\E \jump{E^1(u_h) - u_h} \cdot \avg{\nabla z}}
      \\
      &\leq
      \sup_{g\in \leb{2}(\W), \Norm{g}_{\leb{2}(\W)} \leq 1}
      C\qp{1+h_e^{\qp{\beta-3}/{2}}}\Norm{E^1(u_h) - u_h}_{\leb{2}(\W)} \Norm{\Delta z}_{\leb{2}(\W)}
      \\
      &\leq
      \frac{C\qp{1+h_e^{\qp{\beta-3}/{2}}}}{\gamma*} \Norm{E^1(u_h) - u_h}_{\leb{2}(\W)}
      \\
      &\leq
      C h_e^{1/2} \Norm{\jump{u_h}}_{\leb{2}(\E)},
    \end{split}
  \end{equation}
  for $\beta \geq 3$, concluding the proof.
\end{Proof}

\subsection{Quadrature approximation}

Perhaps the first topic of thought when considering notions of
inconsistency in finite element analysis is that of quadrature
approximations. In the case when the bilinear form contains a positive
definite, non constant diffusion tensor inconsistencies can arise from the
inexact integration \cite[\S 4.1]{Ciarlet:1978}. Take $\A\in\leb{\infty}(\W)^{d\times d}$ and
\begin{equation}
  \label{eq:diff}
  \bi{u}{v} = \int_\W \qp{\A \nabla u} \cdot \nabla v,
\end{equation}
for example, then both linear and bilinear forms are practically
approximated by using quadratures.

\begin{Defn}[Quadrature]
  We introduce 
  \begin{equation}
    Q_K^s(f) := \sum_{i_q = 1}^{n_q} w_{i_q} f(x_{i_q})
  \end{equation}
  to be a quadrature where $\{\qp{x_{i_q}, w_{i_q}}\}_{i_q = 1}^{n_q}$ denote
  the positions and weights of the quadrature points over an element $K$. This is said to
  be an order $s$ quadrature if
  \begin{equation}
    \int_K f = Q_K^s(f) \Foreach f\in\poly{s}(K)
  \end{equation}
 and 
  \begin{equation}
    \norm{\int_K f(x) \d x- Q_K^s(f)} \leq C h^{s+1} \norm{f}_{\sobh{s}(K)} \Foreach f\in\sobh{s}(K).
  \end{equation}
\end{Defn}

Let us, for simplicity, temporarily consider the $\cont{0}$ conforming
method under quadrature approximation applied to (\ref{eq:diff}) which
reads: Find $u_h\in\fes \cap \hoz(\W)$ such that
\begin{equation}
  \label{eq:c0fe}
  \sum_{K\in\T{}} 
  Q_K^s\qp{\qp{\A \nabla u_h} \cdot \nabla v_h}
  = 
  \sum_{K\in\T{}} 
  Q_K^s \qp{f v_h} \Foreach v_h\in\fes \cap \hoz(\W).
\end{equation}

\begin{Rem}[Standard quadrature degree choice]
  A relatively standard choice of quadrature degree is $s=2k-2$. This
  is because if $\A$ is constant then the quadrature allows for exact
  evaluation of the bilinear form. The linear form however is not, in
  general, integrated exactly. Indeed, the linear form can be
  interpreted as
  \begin{equation}
    \label{eq:lin-form-quad}
    l_h(v_h) = \sum_{K\in\T{}} Q_K^{2k-2} (f v_h) = \sum_{K\in\T{}} \int_K P_{k-2} f v_h,
  \end{equation}
  where $P_{k-2}$ denotes the $\leb{2}$ projection into
  $\poly{\max\qp{k-2,0}}(\T{})$.
\end{Rem}

\begin{Lem}[A posteriori bound for quadrature approximations]
  \label{lem:apost-quad}
  Let $u\in\hoz(\W)$ be a weak solution to the variational problem 
  \begin{equation}
    \bi{u}{v} = l(v) \Foreach v\in\hoz(\W),
  \end{equation}
  with $\bi{u}{v}$ given by (\ref{eq:diff}) and $u_h$ solve
  (\ref{eq:c0fe}) with quadrature degree $2k-2$, then
  \begin{equation}
    \Norm{u - u_h}_{\sobh{1}(\W)}
    \leq
    C \qp{\sum_{K\in\T{}} \eta_{R}^2 + \eta_{I}^2 + \sum_{e\in\E} \eta_{J}^2}^{1/2},
  \end{equation}
  where
  \begin{equation}
    \begin{split}
      \eta_R^2 &:= h^2_K \Norm{P_{k-2} f + \div\qp{P_{k-1} \qp{\A \nabla u_h}}}_{\leb{2}(K)}^2
      \\
      \eta_I^2 &:= h^2_K \Norm{f - P_{k-2} f}_{\leb{2}(K)}^2 + \Norm{\qp{P_{k-1} - \Id} \qp{\A \nabla u_h}}_{\leb{2}(K)}^2
      \\
      \eta_J^2 &:= h_e \Norm{\jump{P_{k-1}\qp{\A \nabla u_h}}}_{\leb{2}(e)}^2.
    \end{split}
  \end{equation}
\end{Lem}

\begin{Rem}[Evaluation of the estimator and data oscillation]
  The indicator components $\eta_R, \eta_J$ have the standard form of
  an a posteriori estimator with the exception that they involve
  \emph{projections} of the problem data. This allows for these terms to be
  evaluated \emph{exactly} using the same quadrature as used in the
  scheme, that is, there is no additional approximation involved in
  the computation of these terms.

  The estimator $\eta_I$, representing the inconsistency, is a data
  oscillation term. In general, the integration required in
  computation of the estimator \emph{cannot} be performed exactly
  using the quadrature degree in the scheme.  If some smoothness is
  assumed on the problem data these terms can be considered of higher
  order in comparison to the residual, however in practical
  computations resolution of the problem data is extremely 
  important especially at coarse mesh scales, see \cite{BDN13}. In
  the design of convergent adaptive algorithms the inclusion of these
  terms in the marking strategy is of paramount importance see
  \cite{MorinNochettoSiebert:2002,MekchayNochetto:2005} for further
  arguments.
\end{Rem}

\begin{Proof}[of Lemma \ref{lem:apost-quad}]
  Since the method we consider is conforming Theorem
  \ref{the:the-aposteriori-strang} simplifies to give 
  \begin{equation}
    \begin{split}
      \Norm{\nabla u - \nabla u_h}_{\leb{2}(\W)} 
      &\leq
      \frac{1}{\gamma} 
      \bigg( \sup_{v\in \hoz, \Norm{v}_{\hoz} \leq 1} 
      \qb{
        l_h(v-v_h) - \bih{u_h}{v - v_h} }
      +
      \sup_{v\in\hoz, \Norm{v}_{\hoz} \leq 1} 
      \qb{ l(v) - l_h(v) }
      \\
      &\qquad \qquad
      +
      \sup_{v\in \hoz, \Norm{v}_{\hoz} \leq 1} 
      \qb{
        \bih{u_h}{v} - \bi{u_h}{v} }
      \bigg)
      \\
      & =: 
      \frac{1}{\gamma} \qp{ \E_1 + \E_2 + \E_3 }.
    \end{split}
  \end{equation}
  The first term is a standard residual term and can be controlled
  using standard a posteriori techniques. Indeed, using (\ref{eq:lin-form-quad}) we interpret 
  \begin{equation}
    l_h(v - v_h)
    =
    \sum_{K\in\T{}}
    \int_K
    P_{k-2} f \qp{v - v_h}.
  \end{equation}
  Similarly,
  \begin{equation}
    \bih{u_h}{v-v_h}
    =
    \sum_{K\in\T{}}
    \int_K
    P_{k-1} \qp{\A \nabla u_h} \nabla\qp{v - v_h}.
  \end{equation}
  Hence
  \begin{equation}
    \begin{split}
      l_h(v - v_h)
      -
      \bih{u_h}{v-v_h}
      &=
      \sum_{K\in\T{}}
      \bigg[
      \int_K
      \qp{P_{k-2} f + \div\qp{P_{k-1}\qp{\A\nabla u_h}}} \qp{v-v_h}
      \\
      & \qquad
      - \int_{\partial K}
      P_{k-1}\qp{\A\nabla u_h} \cdot \vec n \qp{v-v_h}
      \bigg]
      \\
      &=
      \sum_{K\in\T{}}
      \bigg[
      \int_K
      \qp{P_{k-2} f + \div\qp{P_{k-1}\qp{\A\nabla u_h}}} \qp{v-v_h}
      \bigg]
      \\
      & \qquad
      - \int_{\E}
      \jump{P_{k-1}\qp{\A\nabla u_h}} \qp{v-v_h}.
    \end{split}
  \end{equation}
  Taking $v_h$ to be the Scott-Zhang interpolant of $v$ and using
  the approximability properties of Proposition \ref{pro:SZ} we have
  \begin{equation}
    \label{eq:E1}
    \begin{split}
      \cE_1
      &\leq
      C \sum_{K\in\T{}}
      \bigg( h_K \Norm{P_{k-2} f + \div\qp{P_{k-1}
          \qp{\A \nabla u_h}}}_{\leb{2}(K)}
      + \sum_{e\in K} h_e^{1/2} \Norm{\jump{P_{k-1}\qp{\A \nabla
            u_h}}}_{\leb{2}(e)} \bigg).
    \end{split}
  \end{equation}
  As for the other term, note that
  \begin{equation}
    \label{eq:E2}
    \begin{split}
      \cE_2
      &= 
      \sup_{v\in\hoz, \Norm{v}_{\hoz} \leq 1} 
      \qb{ l(v) - l_h(v) }
      \\
      &=
      \sup_{v\in\hoz, \Norm{v}_{\hoz} \leq 1} 
      \int_\W \qp{f - P_{k-2} f} v
      \\
      &=
      \sup_{v\in\hoz, \Norm{v}_{\hoz} \leq 1} 
      \int_\W \qp{f - P_{k-2} f} \qp{v - P_{k-2} v}
      \\
      &\leq
      \sup_{v\in\hoz, \Norm{v}_{\hoz} \leq 1} 
      \sum_{K\in\T{}}\Norm{f - P_{k-2} f}_{\leb{2}(K)} \Norm{v - P_{k-2} v}_{\leb{2}(K)}
      \\
      &\leq
      C \sum_{K\in\T{}} h_K \Norm{f - P_{k-2} f}_{\leb{2}(K)},
    \end{split}
  \end{equation}
  using the approximability of the $\leb{2}$ projection operator from
  Proposition \ref{pro:approx-l2proj}.

  The final bound follows from definition since
  \begin{equation}
    \label{eq:E3}
    \begin{split}
      \cE_3
      &=
      \sup_{v\in \hoz, \Norm{v}_{\hoz} \leq 1} 
      \qb{
        \bih{u_h}{v} - \bi{u_h}{v} }
      \\
      &= 
      \sup_{v\in \hoz, \Norm{v}_{\hoz} \leq 1} 
      \sum_{K\in\T{}}
      \int_K
      \qp{P_{k-1} - \Id} \qp{\A \nabla u_h}  \nabla v.
      \\
      &\leq
      \sum_{K\in\T{}}
      \Norm{\qp{P_{k-1} - \Id} \qp{\A \nabla u_h}}_{\leb{2}(K)},
  \end{split}
  \end{equation}
  by Cauchy-Schwarz. Collecting (\ref{eq:E1}),~(\ref{eq:E2}) and
  (\ref{eq:E3}) yields the desired result.
\end{Proof}

\begin{Pro}[A posteriori lower bound under quadrature approximation]
  \label{pro:apost-lower-quad}
  Let $u$, $u_h$ and $\eta_R, \eta_J$ be as in Lemma \ref{lem:apost-quad}  then
  \begin{equation}
    \eta_R + \sum_{e\in K} \eta_J \leq C\qp{\Norm{\nabla u - \nabla u_h}_{\leb{2}(\patch{K})}
      +
      \eta_I
    },
  \end{equation}
  where $\patch{K}$ denotes the set of all elements sharing a common
  edge with $K$.
\end{Pro}
\begin{Proof}
  The proof of this fact is relatively standard, we will include the
  first part, that of the element residual, for completeness and to
  compare with a result in a later section. Throughout this proof we
  will use the convention that $a \lesssim b$ means $a \leq C b$ where
  $C$ is a generic constant that may depend on the problem data, but
  is independent of the meshsize $u$ and $u_h$. We begin by defining
  the bubble functions. Let $b_K$ be the interior bubble function
  defined on the reference triangle $K_{\text{ref}}$ with barycentric
  coordinates $\{\lambda_i\}_{i=0}^d$ through $b_{K_{\text{ref}}} =
  \qp{d+1}^{d+1} \prod_{i = 0}^d \lambda_i$. Due to
  \cite{Verfurth:1996} these functions satisfy
  \begin{equation}
    \label{eq:bubble}
    \Norm{b_K v}_{\leb{2}(K)} 
    \leq
    \Norm{v}_{\leb{2}(K)} 
    \lesssim
    \Norm{b_K v}_{\leb{2}(K)} \text{ for } v\in\poly{k}(K).
  \end{equation}
  In view of this we have that
  \begin{equation}
    \begin{split}
      \Norm{P_{k-2} f + \div\qp{P_{k-1} \qp{\A\nabla u_h}}}^2_{\leb{2}(K)}
      \lesssim
      \Norm{b_K^{1/2} \qp{P_{k-2} f + \div\qp{P_{k-1} \qp{\A\nabla u_h}}}}^2_{\leb{2}(K)},
    \end{split}
  \end{equation}
  so
  \begin{equation}
    \begin{split}
      \Norm{P_{k-2} f + \div\qp{P_{k-1} \qp{\A\nabla u_h}}}_{\leb{2}(K)}^2
      &\lesssim
      \int_K
      \qp{P_{k-2} f - f} b_K  \qp{P_{k-2}f + \div\qp{P_{k-1} \qp{\A\nabla u_h}}}
      \\
      &\qquad +
      \int_K
      \qp{f + \div\qp{{\A\nabla u_h}}} b_K  \qp{P_{k-2} f + \div\qp{P_{k-1} \qp{\A\nabla u_h}}}
      \\
      &\qquad +
      \int_K
      \qp{\div\qp{P_{k-1} \qp{\A\nabla u_h}} - \div\qp{{\A\nabla u_h}}} b_K  \qp{P_{k-2} f + \div\qp{P_{k-1} \qp{\A\nabla u_h}}}
      \\
      &\lesssim
      \int_K
      \qp{P_{k-2} f - f} b_K  \qp{P_{k-2}f + \div\qp{P_{k-1} \qp{\A\nabla u_h}}}
      \\
      &\qquad +
      \int_K
      \qp{\A\nabla u_h - {\A\nabla u}} \cdot \nabla \qp{b_K  \qp{P_{k-2} f + \div\qp{P_{k-1} \qp{\A\nabla u_h}}}}
      \\
      &\qquad +
      \int_K
      \qp{{P_{k-1} \qp{\A\nabla u_h}} - {\A\nabla u_h}} \cdot \nabla \qp{b_K  \qp{P_{k-2} f + \div\qp{P_{k-1} \qp{\A\nabla u_h}}}},
    \end{split}
  \end{equation}
  as $b_K$ vanishes on $\partial K$. Now using Cauchy-Schwarz, further properties of $b_K$ and inverse inequalities
  \begin{equation}
    \begin{split}
      \Norm{P_{k-2} f + \div\qp{P_{k-1} \qp{\A\nabla u_h}}}_{\leb{2}(K)}
      &\lesssim  
      \Norm{
        P_{k-2} f - f
      }_{\leb{2}(K)}
      +
      h_K^{-1}\Norm{\A\nabla u - {\A\nabla u_h}}_{\leb{2}(K)}
      \\
      &\qquad +
      h_K^{-1}\Norm{{{\A\nabla u_h}} - {P_{k-1} \qp{\A\nabla u_h}}}_{\leb{2}(K)}
      \\
      &\lesssim
      \Norm{
        P_{k-2} f - f
      }_{\leb{2}(K)}
      +
      h_K^{-1}\Norm{\nabla u - {\nabla u_h}}_{\leb{2}(K)}
      \\
      &\qquad +
      h_K^{-1}\Norm{{{\A\nabla u_h}} - {P_{k-1} \qp{\A\nabla u_h}}}_{\leb{2}(K)}.
    \end{split}
  \end{equation}
  To see that the jump residual behaves in a similar way one may follow the arguments in \cite[Proof of Thm 10.10]{ErnGuermond:2004} with a slight modification to take into account the data approximation. 
\end{Proof}

The extension to a nonconforming method is almost immediate using the
reconstruction operator from \cite{KarakashianPascal:2003} and the
arguments previously. For instance, 
consider the problem to seek $u_h\in\fes$ such that
\begin{equation}
  \label{eq:dg-quad}
  \bih{u_h}{v_h} = l_h(v_h) \Foreach v_h \in \fes,
\end{equation}
where
\begin{equation}
  \label{eq:bilinear-dg-quad}
  \begin{split}
    \bih{u_h}{v_h}
    &:=
    \sum_{K\in\T{}}
    Q_K^{2k-2}\qp{\qp{\A \nabla u_h} \cdot \nabla v_h}
    -
    \sum_{e\in\E}
    \bigg[
    Q_e^{2k-1}\qp{\jump{u_h}\cdot{\avg{\A\nabla v_h}}}
    \\
    &\qquad\qquad\qquad\qquad\qquad\qquad\qquad +
    Q_e^{2k-1}\qp{\jump{v_h}\cdot{\avg{\A\nabla u_h}}}
    -
    \frac{\sigma}{h} Q_e^{2k}\qp{\jump{v_h}\cdot\jump{u_h}}
    \bigg]
  \end{split}
\end{equation}
and $l_h(v_h)$ is given by (\ref{eq:lin-form-quad}). This is an
example of an IP method after quadrature applied to
(\ref{eq:diff}). Notice the difference in quadrature degree chosen to
integrate over edges. In addition note the similarity in structure to
that of (\ref{eq:IP-mod}). 

\begin{Cor}[A posteriori bound for nonconforming schemes including
  quadrature approximation]
  \label{cor:quad}
  Let $u$ be a weak solution of (\ref{eq:diff}) and $u_h$ solve
  (\ref{eq:dg-quad}), then
  \begin{equation}
    \enorm{u - u_h}{2}
    \leq
    \qp{
    \sum_{K\in\T{}}
      \eta_R^2
      +
      \widetilde{\eta_I}^2
    +
    \sum_{e\in\E}
    \eta_J^2
  }^{1/2}
  \end{equation}
  where
  \begin{equation}
    \begin{split}
      \eta_R^2 &:= h^2_K \Norm{P_{k-2} f + \div\qp{P_{k-1} \qp{\A \nabla u_h}}}_{\leb{2}(K)}^2
      \\
      \widetilde{\eta_I}^2 &:= h^2_K \Norm{f - P_{k-2} f}_{\leb{2}(K)}^2 + \Norm{\qp{P_{k-1} - \Id} \qp{\A \nabla E^1(u_h)}}_{\leb{2}(K)}^2
      \\
      \eta_J^2 &:= h_e \Norm{\jump{P_{k-1}\qp{\A \nabla u_h}}}_{\leb{2}(e)}^2 + \sigma h_e^{-1} \Norm{\jump{u_h}}_{\leb{2}(e)}^2.
    \end{split}
  \end{equation}
  Also for $\sigma$ large enough
  \begin{equation}
    \eta_R + \sum_{e\in K} \eta_J
    \leq
    C \qp{
    \Norm{\nabla u - \nabla u_h}_{\leb{2}(K)} +  {\eta_I}}
  \end{equation}
  with
  \begin{equation}
    \eta_I^2 := h^2_K \Norm{f - P_{k-2} f}_{\leb{2}(K)}^2 + \Norm{\qp{P_{k-1} - \Id} \qp{\A \nabla u_h}}_{\leb{2}(K)}^2
  \end{equation}
\end{Cor}
\begin{Proof}
  The main difference between control of the conforming approximation
  and the nonconforming is that of the penalty terms. To see these
  control the error from below we modify the arguement given in
  \cite[Lem 5.30]{Di-PietroErn:2012}. The main idea is that for
  variational second order problems there is a conformal
  reconstruction in the dG space. This allows us to make use of the
  numerical scheme in the following fashion. We interpret (\ref{eq:bilinear-dg-quad}) as
  \begin{equation}
    \begin{split}
      \bih{u_h}{v_h} &= \int_\W P_{k-1} \qp{\A \nabla_h u_h} \cdot \nabla_h v_h
      \\&\qquad -
      \int_\E 
      \bigg[
      \jump{u_h} \avg{P_{k-1} \qp{\A\nabla v_h}}
      +
      \jump{v_h} \avg{P_{k-1} \qp{\A\nabla u_h}}
      -
      \sigma h_e^{-1} \jump{u_h} \cdot \jump{v_h}
    \bigg],
    \end{split}
  \end{equation}
  which, using the fact that $u_h$ solves (\ref{eq:dg-quad}), allows
  us to choose $v_h = u_h - E^1(u_h)$. From this we see, since
  $\jump{E^1(u_h)} = 0$, that after an integration by parts
  \begin{equation}
    \label{eq:jump-proof-0}
    \begin{split}
      \sigma h_e^{-1} \Norm{\jump{u_h}}^2_{\leb{2}(\E)}
      &=
      \int_\W \qp{P_{k-2} f + \div\qp{P_{k-1}\qp{\A\nabla u_h}}} \qp{u_h - E^1(u_h)}
      \\
      &\qquad 
      - \int_\E 
      \bigg[
      \jump{P_{k-1} \qp{\A\nabla u_h}} \cdot \avg{u_h - E^1(u_h)}
      -
      \jump{u_h}\cdot \avg{P_{k-1} \qp{\A \nabla \qp{u_h - E^1(u_h)}}}
      \bigg]
      \\
      &\leq
      \Norm{{P_{k-2} f + \div\qp{P_{k-1}\qp{\A\nabla u_h}}}}_{\leb{2}(\W)}
      \Norm{u_h - E^1(u_h)}_{\leb{2}(\W)}
      \\
      &\qquad +
      \Norm{\jump{P_{k-1} \qp{\A\nabla u_h}}}_{\leb{2}(\E)} \Norm{\avg{u_h - E^1(u_h)}}_{\leb{2}(\E)}
      \\
      &\qquad +
      \Norm{\jump{u_h}}_{\leb{2}(\E)} \Norm{\avg{P_{k-1} \qp{\A \nabla \qp{u_h - E^1(u_h)}}}}_{\leb{2}(\E)}
      \\
      &=: \cJ_1 + \cJ_2 + \cJ_3.
    \end{split}
  \end{equation}
  The first term, owing to the properties of the reconstruction
  $E^1(u_h)$ is controlled by
  \begin{equation}
    \label{eq:jump-proof-1}
    \begin{split}
      \cJ_1 &=
      \Norm{{P_{k-2} f + \div\qp{P_{k-1}\qp{\A\nabla u_h}}}}_{\leb{2}(\W)}
      \Norm{u_h - E^1(u_h)}_{\leb{2}(\W)}
      \\
      &\leq
      Ch_e^{1/2}
      \Norm{{P_{k-2} f + \div\qp{P_{k-1}\qp{\A\nabla u_h}}}}_{\leb{2}(\W)}
      \Norm{\jump{u_h}}_{\leb{2}(\E)}.
    \end{split}
  \end{equation}
  The second term, via a trace trace inequality
  \begin{equation}
    \label{eq:jump-proof-2}
    \begin{split}
      \cJ_2 &=
      \Norm{\jump{P_{k-1} \qp{\A\nabla u_h}}}_{\leb{2}(\E)} \Norm{\avg{u_h - E^1(u_h)}}_{\leb{2}(\E)}
      \\
      &\leq
      C h_e^{-1/2}\Norm{\jump{P_{k-1} \qp{\A\nabla u_h}}}_{\leb{2}(\E)} \Norm{{u_h - E^1(u_h)}}_{\leb{2}(\W)}
      \\
      &\leq
      C \Norm{\jump{P_{k-1} \qp{\A\nabla u_h}}}_{\leb{2}(\E)} \Norm{\jump{u_h}}_{\leb{2}(\E)}.
    \end{split}
  \end{equation}
  Finally the third term, again by a trace inequality and the
  stability of the $\leb{2}$ projection
  \begin{equation}
    \label{eq:jump-proof-3}
    \begin{split}
      \cJ_3
      &=
      \Norm{\jump{u_h}}_{\leb{2}(\E)} \Norm{\avg{P_{k-1} \qp{\A \nabla \qp{u_h - E^1(u_h)}}}}_{\leb{2}(\E)}
      \\
      &\leq
      C h_e^{-1/2}
      \Norm{\jump{u_h}}_{\leb{2}(\E)} \Norm{{P_{k-1} \qp{\A \nabla \qp{u_h - E^1(u_h)}}}}_{\leb{2}(\W)}
      \\
      &\leq
      C h_e^{-1/2}
      \Norm{\jump{u_h}}_{\leb{2}(\E)} \Norm{\A \nabla \qp{u_h - E^1(u_h)}}_{\leb{2}(\W)}
      \\
      &\leq
      \widetilde{C} h_e^{-1}
      \Norm{\jump{u_h}}_{\leb{2}(\E)}^2.
   \end{split}
  \end{equation}
  Substituting (\ref{eq:jump-proof-1}), (\ref{eq:jump-proof-2}) and
  (\ref{eq:jump-proof-3}) into (\ref{eq:jump-proof-0}) yields
  \begin{equation}
    \begin{split}
      \qp{\sigma - \widetilde C} h_e^{-1} \Norm{\jump{u_h}}^2_{\leb{2}(\E)}
      &\leq
      C \qp{\Norm{{P_{k-2} f + \div\qp{P_{k-1}\qp{\A\nabla u_h}}}}_{\leb{2}(\W)}
      +
      \Norm{\jump{P_{k-1} \qp{\A\nabla u_h}}}_{\leb{2}(\E)}}.
    \end{split}
  \end{equation}
  Hence for $\sigma > \widetilde C$ the penalty term is controlled by
  the other residual terms which in turn, using the arguements from
  Proposition \ref{pro:apost-lower-quad}, are controlled by the error
  up to the data oscillation terms given by $\eta_I$, concluding the
  proof.
\end{Proof}

\begin{Rem}[Inconsistencies require computation of reconstructions]
  In the statement of the bound for the nonconforming approximation of
  the Laplacian in Lemma \ref{lem:lap-apost} there is no requirement
  to compute the reconstruction operator, $E^1$. Indeed, if there is no
  inconsistency in the approximation of the nonconstant diffusion
  problem (\ref{eq:diff}) then it is still not necessary to actually
  compute the reconstruction operator. The nature of the inconsistency
  term $\widetilde{\eta_I}$ appearing in Corollary \ref{cor:quad} requires
  computation of the reconstruction of $u_h$. We will see in a later
  section that inconsistencies do not \emph{always} require
  computation of this operator.
\end{Rem}

\section{Unbalanced variational and nonvariational problems}
\label{sec:unbalanced}
In this section we will take a different route. For simplicity
initially we shall ignore any possible effect of quadrature and focus
on the effects of inconsistencies in unbalanced problems. We will
consider the case $W = \sobh{2}(\W) \cap
\hoz(\W)$ and $V = \leb{2}(\W)$. Let
\begin{equation}
  \label{eq:nondiv}
  \bi{u}{v} = \int_\W \frob{\A}{\Hess u},
\end{equation}
where $\A \in \reals^{d\times d}$ is positive definite diffusion
tensor, $\Hess u$ denotes the Hessian matrix of $u$ and $\frob{\geomat
  X}{\geomat Y} = \trace\qp{\Transpose{\geomat X}\geomat Y}$ is the
Frobenious product between matrices. The operator (\ref{eq:nondiv}) is
in nondivergence form. Galerkin schemes falling into the framework
presented here are given in
\cite{LakkisPryer:2011a,DednerPryer:2013}. Both are nonconforming,
inconsistent approximation schemes. 
\begin{Hyp}[Strong solutions exist]
  \label{hyp:strong}
  Under either smoothness assumptions on $\A$ \cite{Gilbarg:1983},
  Campanato \cite{Campanato:1994} or Cordes \cite{Cordes:1961}
  conditions on $\A$ the problem
  \begin{equation}
    \label{eq:non-div-eqn}
    \bi{u}{v} = \int_\W f v
  \end{equation}
  has a strong solution, that is, there exists a constant $C>0$ such that
  \begin{equation}
    \Norm{u}_{\sobh{2}(\W)} \leq C \Norm{f}_{\leb{2}(\W)}.
  \end{equation}
  This is equivalent to Assumption \ref{ass:cont-inf-sup}. Note that
  if $\A = -\geomat I$ then this is a natural framework to study the
  Laplacian if it has strong solutions.
\end{Hyp}

\subsection{Consistent nonvariational methods}

Henceforth we will take $k \geq 2$. For
$\A\in\leb{\infty}(\W)^{d\times d}$ satisfying Assumption
\ref{hyp:strong} we define
\begin{equation}
  \label{eq:cons-nonvar-bilinear}
  \bih{u_h}{v_h} 
  =
  \int_\W \frob{\A}{\Hess_h u_h} \ v_h
  -
  \int_\E \qp{\frob{\tjump{\nabla u_h}}{\avg{ \A v_h}} 
    -
    \sigma h_e \jump{\nabla u_h}\jump{\nabla v_h}
    -
    \sigma h_e^{-1} \jump{u_h}\cdot\jump{v_h}
  },
\end{equation}
and seek $u_h\in\fes$ such that
\begin{equation}
  \label{eq:non-div-cons}
  \bih{u_h}{v_h} = \int_\W f v_h \Foreach v_h \in \fes.
\end{equation}
Note for $u\in\sobh{2}(\W)$ the method is consistent, in that
\begin{equation}
  \bih{u}{v} = \bi{u}{v} \Foreach v\in\sobh{2}(\T{}).
\end{equation}

\begin{Rem}[Extension of the bilinear form]
  \label{rem:extension-nonvar}
  The discrete bilinear form (\ref{eq:cons-nonvar-bilinear}) only
  makes sense over $\sobh2\qp{\T{}}\times \sobh2\qp{\T{}}$. To make
  use of the a posteriori framework in \S\ref{sec:setup} we require an
  extension to ensure the appropriate stability arguemnts can be applied. That
  is we require the bilinear form to be extended to
  $\sobh2\qp{\T{}}\times \leb 2\qp{\W}$. To do this for
  $\qp{u,v}\in\sobh2\qp{\T{}}\times \leb 2\qp{\W}$ we define
  \begin{equation}
    \label{eq:cons-nonvar-bilinear-mod}
    \bih{u}{v} 
    =
    \int_\W \frob{\A}{\Hess_h u} \ v
    -
    \int_\E \qp{\frob{\tjump{\nabla u}}{\avg{ \A P_k v}} 
      -
      \sigma h_e \jump{\nabla u}\jump{\nabla\qp{ P_k v}}
      -
      \sigma h_e^{-1} \jump{u}\cdot\jump{P_k v}
    }.
  \end{equation}
  Notice that the modified bilinear form
  (\ref{eq:cons-nonvar-bilinear-mod}) coincides with
  (\ref{eq:cons-nonvar-bilinear}) over $\fes\times\fes$ and that it
  satisfies
  \begin{equation}
    \bih{u}{v} \leq C \eenorm{u}{2} \Norm{v}_{\leb{2}(\W)} \text{ for } \qp{u,v}\in\qp{\sobh2\qp{\T{}}\times\leb{2}(\W)}.
  \end{equation}
\end{Rem}

\begin{The}[Primal consistent a posteriori upper bound]
  \label{eq:primal-cons-upper}
  Let $u\in\sobh{2}(\W)$ solve (\ref{eq:non-div-eqn}) and $u_h\in\fes$
  be the finite element approximation given by (\ref{eq:non-div-cons}), then
  \begin{equation}
    \eenorm{u - u_h}{2} 
    \leq 
    C 
    \qp{\sum_{K\in\T{}}\eta_R^2 + \sum_{e\in\E}\eta_J^2}^{1/2},
  \end{equation}
  where
  \begin{gather}
    \eta_R^2 := \Norm{f - \frob{\A}{\Hess u_h}}_{\leb{2}(K)}^2
    \\
    \eta_J^2 := 
    h_e^{-1} \Norm{\jump{\nabla u_h}}_{\leb{2}(e)}^2
    +
    h_e^{-3} \Norm{\jump{u_h}}_{\leb{2}(e)}^2.
  \end{gather}
\end{The}
\begin{Proof}
  Making use of the modified bilinear form given in Remark
  \ref{rem:extension-nonvar} we may apply the framework of Theorem
  \ref{the:the-aposteriori-strang} directly here and we have
  \begin{equation}
    \label{eq:aposteriori-primal}
    \begin{split}
      \Norm{e}_{\sobh2(\W)}
      &\leq
      \frac{1}{\gamma}
      \bigg(
      \sup_{v\in {\leb{2}(\W)}, \Norm{v}_{\leb{2}(\W)} \leq 1}
      \qb{
        l_h(v-v_h) - \bih{u_h}{v - v_h}
      }
      +
      \sup_{v\in {\leb{2}(\W)}, \Norm{v}_{\leb{2}(\W)} \leq 1}
      \qb{
        l(v) - l_h(v)
      }
      \\ 
      &\qquad \qquad +
      \sup_{v\in {\leb{2}(\W)}, \Norm{v}_{\leb{2}(\W)} \leq 1}
      \qb{
        \bih{E^2(u_h)}{v} - \bi{E^2(u_h)}{v}
      }
      \bigg)
      +
      \qp{1+ \frac{C_B}{\gamma}} \Norm{e^N}_{\sobh2(\W)}
      \\
      & =: 
      \frac{1}{\gamma} \qp{ \E_1 + \E_2 + \E_3 + \qp{\gamma + C_B }\E_4}.
    \end{split}
  \end{equation}
  Note that $\cE_2 = 0$ and since $E^2(u_h)\in\sobh{2}(\W)$ all jump
  terms from (\ref{eq:cons-nonvar-bilinear-mod}) vanish and we see
  $\cE_3 = 0$. To control $\cE_1$ we see in view of Cauchy-Schwarz
  \begin{equation}
    \begin{split}
      \cE_1
      &=
      \sup_{v\in {\leb{2}(\W)}, \Norm{v}_{\leb{2}(\W)} \leq 1}
      \qb{
        l_h(v-v_h) - \bih{u_h}{v - v_h}
      }
      \\
      &=
      \sup_{v\in {\leb{2}(\W)}, \Norm{v}_{\leb{2}(\W)} \leq 1}
      \bigg[
      \int_\W \qp{f - \frob{\A}{\Hess_h u_h}} \qp{v - v_h}
      +
      \int_\E \frob{\tjump{\nabla u_h}}{\avg{\A P_k\qp{v-v_h}}}
      \\
      &\qquad\qquad\qquad -
      \sigma h_e\int_\E \jump{\nabla u_h} \jump{\nabla \qp{P_k \qp{v-v_h}}}
      -
      \sigma h_e^{-1}\int_\E \jump{u_h}\cdot\jump{{P_k \qp{v-v_h}}}
      \bigg]
      \\
      &\leq
      \sup_{v\in {\leb{2}(\W)}, \Norm{v}_{\leb{2}(\W)} \leq 1}
      \bigg[
      \sum_{K\in\T{}}
      \Norm{f - \frob{\A}{\Hess_h u_h}}_{\leb{2}(K)}
      \Norm{v - v_h}_{\leb{2}(K)}
      \\
      &\qquad +
      \sum_{e\in\E}
      \bigg(
      \Norm{\jump{\nabla u_h}}_{\leb{2}(e)}
      \qp{
        \Norm{\avg{\A P_k\qp{v - v_h}}}_{\leb{2}(e)}
        +
        \sigma h_e\Norm{\jump{\nabla\qp{P_k\qp{v - v_h}}}}_{\leb{2}(e)}
      }
      \\
      &\qquad \qquad 
      +
      \sigma h_e^{-1}\Norm{\jump{u_h}}_{\leb{2}(e)}
      \Norm{\jump{P_k\qp{v-v_h}}}_{\leb{2}(e)}
      \bigg)\bigg].
    \end{split}
  \end{equation}
  Now choosing $v_h=0$, using trace and inverse inequalities and the
  stability of $P_k$ in $\leb{2}$ we have
  \begin{equation}
    \begin{split}
      \cE_1
      &\leq
      \sup_{v\in {\leb{2}(\W)}, \Norm{v}_{\leb{2}(\W)} \leq 1}
      \sum_{K\in\T{}}
      \Bigg(
      \Norm{f - \frob{\A}{\Hess_h u_h}}_{\leb{2}(K)}
      \Norm{v}_{\leb{2}(K)}
      \\
      &\qquad +
      C
      \sum_{e\in K}
      \bigg(
      \Norm{\jump{\nabla u_h}}_{\leb{2}(e)}
      \qp{
        h_e^{-1/2}\Norm{{P_k\qp{v}}}_{\leb{2}(K)}
        +
        \sigma h_e^{1/2}\Norm{{\nabla{P_k v}}}_{\leb{2}(K)}
      }
      \\
      &\qquad \qquad +
      \Norm{\jump{u_h}}_{\leb{2}(e)}
      h_e^{-3/2}\Norm{{P_k{v}}}_{\leb{2}(K)}
      \bigg)
      \Bigg)
      \\
      &\leq
      \sup_{v\in {\leb{2}(\W)}, \Norm{v}_{\leb{2}(\W)} \leq 1}
      \Norm{v}_{\leb{2}(\W)}
      \sum_{K\in\T{}}
      \Bigg(
      \Norm{f - \frob{\A}{\Hess_h u_h}}_{\leb{2}(K)}
      \\
      &\qquad +
      C
      \sum_{e\in K}
      \bigg(
      C
      h_e^{-1/2}
      \Norm{\jump{\nabla u_h}}_{\leb{2}(e)}
      +
      h_e^{-3/2}\Norm{\jump{u_h}}_{\leb{2}(e)}
      \bigg)
      \Bigg)
      \\
      &\leq
      \sum_{K\in\T{}}
      \Bigg(
      \Norm{f - \frob{\A}{\Hess_h u_h}}_{\leb{2}(K)}
      +
      C
      \sum_{e\in K}
      \bigg(
      C
      h_e^{-1/2}
      \Norm{\jump{\nabla u_h}}_{\leb{2}(e)}
      +
      h_e^{-3/2}\Norm{\jump{u_h}}_{\leb{2}(e)}
      \bigg)
      \Bigg).
    \end{split}
  \end{equation}
  To conclude we use the reconstruction bounds given in Lemma
  \ref{lem:reconstruction-bounds} to control $\cE_4$ as required.
\end{Proof}

\begin{Pro}[Lower bounds]
  \label{pro:cons-lower-bounds}
  Under the conditions of Theorem \ref{eq:primal-cons-upper} the
  following global a posteriori lower bound holds:
  \begin{equation}
    \sum_{K\in\T{}}\eta_R + \sum_{e\in\E}\eta_J
    \leq
    \Norm{\Hess u - \Hess_h u_h}_{\leb{2}(\W)},
  \end{equation}
  where
  \begin{gather}
    \eta_R^2 := \Norm{f - \frob{\A}{\Hess u_h}}_{\leb{2}(K)}^2
    \\
    \eta_J^2 := 
    h_e^{-1} \Norm{\jump{\nabla u_h}}_{\leb{2}(e)}^2
    +
    h_e^{-3} \Norm{\jump{u_h}}_{\leb{2}(e)}^2.
  \end{gather}
\end{Pro}
\begin{Proof}
  The bound on the interior residual is trivial since the error and
  residual are being measured in the same norm there is no need to use
  bubble functions to invoke inverse inequalities. Indeed,
  \begin{equation}
    \label{eq:pflb1}
    \Norm{f - \frob{\A}{\Hess_h u_h}}_{\leb{2}(\W)}
    \leq
    \Norm{A}_{\leb{\infty}(\W)} \Norm{\Hess u - \Hess_h u_h}_{\leb{2}(\W)}.
  \end{equation}
  For the jump terms we may use a similar argument to that given in
  the proof of Corollary \ref{cor:quad}. We choose $v_h = u_h -
  E^1(u_h)$ in (\ref{eq:cons-nonvar-bilinear}) and then by definition of the scheme
  \begin{equation}
    \begin{split}
      \sigma h_e^{-1}
      \Norm{\jump{u_h}}_{\leb{2}(e)}^2
      &\leq
      \sigma h_e^{-1}
      \Norm{\jump{\nabla u_h}}_{\leb{2}(\E)}
      \Norm{\jump{\nabla u_h - \nabla E^1(u_h)}}_{\leb{2}(\E)}
      +
      \Norm{\jump{\nabla u_h}}_{\leb{2}(\E)}
      \Norm{\avg{u_h - E^1(u_h)}}_{\leb{2}(\E)}
      \\
      &\qquad +
      \Norm{f - \frob{\A}{\Hess_h u_h}}_{\leb{2}(\W)}
      \Norm{{u_h - E^1(u_h)}}_{\leb{2}(\W)}.
    \end{split}
  \end{equation}
  Now using the properties of $E^1$ from Lemma
  \ref{lem:reconstruction-bounds} we have
  \begin{equation}
    \begin{split}
      \sigma h_e^{-1}
      \Norm{\jump{u_h}}_{\leb{2}(\E)}^2
      &\leq
      C \sigma
      \Norm{\jump{\nabla u_h}}_{\leb{2}(\E)}
      \Norm{\jump{u_h}}_{\leb{2}(\E)}
      +
      h_e^{1/2}
      \Norm{f - \frob{\A}{\Hess u_h}}_{\leb{2}(K)}
      \Norm{\jump{u_h}}_{\leb{2}(\E)}.
    \end{split}
  \end{equation}
  Now using the fact that
  \begin{equation}
    \label{eq:funny-jump-bound}
    h_e^{-1} \Norm{\jump{\nabla u_h}}_{\leb{2}(\E)}^2
    =
    h_e^{-1} \Norm{\jump{\nabla u - \nabla u_h}}_{\leb{2}(\E)}^2
    \leq
    \Norm{\Hess u - \Hess_h u_h}_{\leb{2}(\W)}^2
  \end{equation}
  together with (\ref{eq:pflb1}) yields
  \begin{equation}
    h_e^{-3/2}
    \Norm{\jump{u_h}}_{\leb{2}(\E)}
    +
    h_e^{-1/2} \Norm{\jump{\nabla u_h}}_{\leb{2}(\E)}
    \leq
    C \Norm{\Hess u - \Hess_h u_h}_{\leb{2}(\W)},
  \end{equation}
  concluding the proof.
\end{Proof}

\subsection{Inconsistent nonvariational methods}

We now introduce an inconsistent method for the approximation of
(\ref{eq:nondiv}). This is a slight modification of that presented in
\cite{DednerPryer:2013} and related to \cite{LakkisPryer:2011a}. The
method we consider is given as follows: Let $\H(u_h)\in \fes^{d\times
  d}$ be such that for all $\phi\in \fes$
\begin{equation}
  \label{eq:defn-of-H}
  \int_\W \H(u_h) \phi
  =
  \int_\W \Hess_h u_h \phi - \int_\E \tjump{\nabla u_h} \avg{\phi} - \jump{u_h} \otimes \avg{\nabla \phi},
\end{equation}
then we seek the tuple $\qp{u_h, \H(u_h)}$ such that
\begin{equation}
  \label{eq:non-div-fem}
  \bih{u_h}{v_h} = \int_\W f v_h \Foreach v_h \in \fes\cap\hoz(\T{}),
\end{equation}
with
\begin{equation}
  \bih{u_h}{v_h} := \int_\W \frob{\A}{\H(u_h)} v_h + \int_\E \sigma h_e \jump{\nabla u_h} \jump{\nabla v_h} + \int_\E \sigma h_e^{-1} \jump{u_h}\cdot \jump{v_h}.
\end{equation}

\begin{The}[Stability of $\H(\cdot)$ {\cite[Thm 4.11]{Pryer:2014}}]
  \label{lem:stability-of-H}
  Let $\H(\cdot)$ be defined as in (\ref{eq:defn-of-H}), then it is
  stable in the sense that
  \begin{equation}
    \label{eq:bound-the-liftings}
    \begin{split}
      \Norm{\Hess_h v_h - \H(v_h)}^2_{\leb{2}(\W)^{d\times d}} 
      &\leq
      C\qp{
        \int_{\E} h_e^{-1} 
        \norm{\jump{\nabla_h v_h}}^2
        +
        h_e^{-3}
        \norm{\jump{v_h}}^2
      }.
    \end{split}
  \end{equation}
  Consequently we have
  \begin{equation}
    \label{eq:stab-of-H}
    \Norm{\H(v_h)}^2_{\leb{2}(\W)^{d\times d}}  
    \leq
    C \eenorm{v_h}{2}^2.
  \end{equation}
\end{The}

\begin{The}[Primal a posteriori error upper bound]
  \label{the:primal-upper}
  Let $u\in\sobh{2}(\W)$ solve (\ref{eq:non-div-eqn}) and $u_h\in\fes$
  be the finite element approximation given by (\ref{eq:non-div-fem}), then
  \begin{equation}
    \eenorm{u - u_h}{2} 
    \leq 
    C 
    \qp{\sum_{K\in\T{}}\eta_R^2 + \sum_{e\in\E}\eta_J^2}^{1/2},
  \end{equation}
  where
  \begin{gather}
    \eta_R^2 := \Norm{f - \frob{\A}{\H(u_h)}}_{\leb{2}(K)}^2
    \\
    \eta_J^2 := 
    h_e^{-1} \Norm{\jump{\nabla u_h}}_{\leb{2}(e)}^2
    +
    h_e^{-3} \Norm{\jump{u_h}}_{\leb{2}(e)}^2.
  \end{gather}
\end{The}
\begin{Proof}
  The proof is similar to that of Theorem
  \ref{eq:primal-cons-upper}. Using the abstract result of Theorem
  \ref{the:the-aposteriori-strang} we proceed to bound each $\E_i$
  term. In view of Cauchy--Schwartz we have
  \begin{equation}
    \begin{split}
      \E_1 &= \sup_{v\in \leb{2}(\W), \Norm{v}_{\leb{2}(\W)}\leq 1}\int_\W \qp{f - \frob{\A}{\H(u_h)}}\qp{v-v_h}
      \\
      &\leq
      \sup_{v\in \leb{2}(\W), \Norm{v}_{\leb{2}(\W)}\leq 1}
      \sum_{K\in\T{}}
      \Norm{f - \frob{\A}{\H(u_h)}}_{\leb{2}(K)} \Norm{v - v_h}_{\leb{2}(K)}.
    \end{split}
  \end{equation}
  Choosing $v_h = 0$ we have
  \begin{equation}
    \label{eq:primal-nonvar-1}
    \E_1 \leq C
      \sum_{K\in\T{}}
      \Norm{f - \frob{\A}{\H(u_h)}}_{\leb{2}(K)}
      .
  \end{equation}
  The consistency term $\E_2 = 0$. For the second consistency term we
  have that
  \begin{equation}
    \begin{split}
      \E_3 
      &=
      \sup_{v\in \leb{2}(\W), \Norm{v}_{\leb{2}(\W)}\leq 1}
      \frac{1}{\gamma}
      \int_\W \frob{\A}{\qp{\H(E^2(u_h)) - \Hess E^2(u_h)}} v
      \\
      &\leq
      C \Norm{\A}_{\leb{\infty}(\W)} \sum_{K\in\T{}} \Norm{\H(E^2(u_h)) - \Hess E^2(u_h)}_{\leb{2}(K)}    
      \\
      &\leq
      C \Norm{\A}_{\leb{\infty}(\W)} \sum_{K\in\T{}} 
      \bigg(
      \Norm{\H(E^2(u_h)) - \H(u_h)}_{\leb{2}(K)} 
      +
      \Norm{\H(u_h) - \Hess_h u_h}_{\leb{2}(K)}
      \\
      &\qquad \qquad \qquad \qquad  \qquad \qquad  \qquad \qquad +
      \Norm{\Hess_h u_h - \Hess E^2(u_h)}_{\leb{2}(K)} \bigg).
      \end{split}
  \end{equation}
  Making use of Lemmata \ref{lem:stability-of-H} and \ref{lem:reconstruction-bounds}, we see that 
  \begin{equation}
    \begin{split}
      \sum_{K\in\T{}}\Norm{\H(E^2(u_h)) - \H(u_h)}_{\leb{2}(K)}
      &\leq
      C\eenorm{E^2(u_h) - u_h}{2}
      \\
      &\leq
      C
      \qp{
        \sum_{e \in \E}
        h_e^{-1} 
        \Norm{\jump{\nabla_h v_h}}^2_{\leb{2}(e)}
        +
        h_e^{-3}
        \Norm{\jump{v_h}}^2_{\leb{2}(e)}
      }^{1/2}.
    \end{split}
  \end{equation}
  Again from Lemma \ref{lem:stability-of-H} we have that
  \begin{equation}
    \sum_{K\in\T{}}
    \Norm{\H(u_h) - \Hess_h u_h}_{\leb{2}(K)}
    \leq
    C
    \qp{
      \sum_{e \in \E}
      h_e^{-1} 
      \Norm{\jump{\nabla_h v_h}}^2_{\leb{2}(e)}
      +
      h_e^{-3}
      \Norm{\jump{v_h}}^2_{\leb{2}(e)}
    }^{1/2}.
  \end{equation}
  and Lemma \ref{lem:reconstruction-bounds} gives that
  \begin{equation}
    \sum_{K\in\T{}}
    \Norm{\Hess_h u_h - \Hess E^2(u_h)}_{\leb{2}(K)}
    \leq
    C
    \qp{
      \sum_{e \in \E}
      h_e^{-1} 
      \Norm{\jump{\nabla_h v_h}}^2_{\leb{2}(e)}
      +
      h_e^{-3}
      \Norm{\jump{v_h}}^2_{\leb{2}(e)}
    }^{1/2}.
  \end{equation}
  Hence we have that
  \begin{equation}
    \label{eq:primal-nonvar-2}
    \begin{split}
      \E_3 
      &\leq
          C
    \qp{
      \sum_{e \in \E}
      h_e^{-1} 
      \Norm{\jump{\nabla_h v_h}}^2_{\leb{2}(e)}
      +
      h_e^{-3}
      \Norm{\jump{v_h}}^2_{\leb{2}(e)}
    }^{1/2}.
    \end{split}
  \end{equation}
  We conclude by noting that from Lemma \ref{lem:reconstruction-bounds}
  \begin{equation}
    \label{eq:primal-nonvar-3}
    \E_4 = \eenorm{E(u_h) - u_h}{2}
    \leq
    C
    \qp{
      \sum_{e \in \patch{K}}
      h_e^{-1} 
      \Norm{\jump{\nabla_h v_h}}^2_{\leb{2}(e)}
      +
      h_e^{-3}
      \Norm{\jump{v_h}}^2_{\leb{2}(e)}
    }^{1/2}.
  \end{equation}
  Collecting (\ref{eq:primal-nonvar-1}), (\ref{eq:primal-nonvar-2})
  and (\ref{eq:primal-nonvar-3}) yields the desired result.
\end{Proof}

\begin{Pro}[Primal a posteriori lower bound]
  \label{the:primal-lower}
  Using the notation of Theorem \ref{the:primal-upper} we have for
  $\sigma$ large enough a lower bound of the form:
  \begin{equation}
    \sum_{K\in\T{}}\eta_R + \sum_{e\in\E} \eta_J
    \leq 
    C \qp{\Norm{\Hess u - \Hess u_h}_{\leb{2}(\W)}}.
  \end{equation}
\end{Pro}
\begin{Proof}
  Notice, as in the consistent case there is no data oscillation, this
  is since the inconsistency can be linked to the nonconformity via
  Lemma \ref{lem:stability-of-H}. We begin by noting
  \begin{equation}
    \begin{split}
      \Norm{f - \frob{\A}{\H(u_h)}}_{\leb{2}(\W)} 
      &\leq
      \Norm{\frob{\A}{\qp{\Hess u - \Hess_h u_h}}}_{\leb{2}(\W)} 
      +
      \Norm{\frob{\A}{\qp{\Hess_h u_h - \H(u_h)}}}_{\leb{2}(\W)} 
      \\
      &\leq
      C 
      \qp{
        \Norm{{{\Hess u - \Hess_h u_h}}}_{\leb{2}(\W)} 
        +
        \Norm{{{\Hess_h u_h - \H(u_h)}}}_{\leb{2}(\W)} 
      }
      \\
      &\leq
      C 
      \qp{
        \Norm{{{\Hess u - \Hess_h u_h}}}_{\leb{2}(\W)} 
        +
        h_e^{-1/2} 
        \Norm{\jump{\nabla u_h}}_{\leb{2}(\E)}
        +
        h_e^{-3/2}
        \Norm{\jump{u_h}}_{\leb{2}(\E)}        
      },
    \end{split}
  \end{equation}
  in view of Lemma \ref{lem:stability-of-H}. Now, as in the proof of
  Proposition \ref{pro:cons-lower-bounds} making use of the scheme
  (\ref{eq:non-div-fem}) with $v_h = u_h - E^1(u_h)$ we see
  \begin{equation}
    \begin{split}
      \sigma h_e^{-1}
      \Norm{\jump{u_h}}_{\leb{2}(e)}^2
      &\leq
      \sigma h_e^{-1}
      \Norm{\jump{\nabla u_h}}_{\leb{2}(\E)}
      \Norm{\jump{\nabla u_h - \nabla E^1(u_h)}}_{\leb{2}(\E)}
      \\
      &\qquad +
      \Norm{f - \frob{\A}{\H(u_h)}}_{\leb{2}(\W)}
      \Norm{{u_h - E^1(u_h)}}_{\leb{2}(\W)}
      \\
      &\leq
      C
      \Norm{\jump{\nabla u_h}}_{\leb{2}(\E)}
      \Norm{\jump{u_h}}_{\leb{2}(\E)}
      \\&\qquad + 
      Ch^{1/2}
      \Norm{\jump{u_h}}_{\leb{2}(\E)}        
      \qp{
        \Norm{{{\Hess u - \Hess_h u_h}}}_{\leb{2}(\W)} 
        +
        h_e^{-1/2} 
        \Norm{\jump{\nabla u_h}}^2_{\leb{2}(\E)}
        +
        h_e^{-3/2}
        \Norm{\jump{u_h}}^2_{\leb{2}(\E)}        
      }.
    \end{split}
  \end{equation}
  Rearranging the inequality we have 
  \begin{equation}
    \begin{split}
      \qp{\sigma -C} h_e^{-3/2}
      \Norm{\jump{u_h}}_{\leb{2}(e)}
      &\leq
      C
      h_e^{1/2}
      \Norm{\jump{\nabla u_h}}_{\leb{2}(\E)}
      + 
      C
      \Norm{{{\Hess u - \Hess_h u_h}}}_{\leb{2}(\W)} 
      .
    \end{split}
  \end{equation}
  Choosing $\sigma$ large enough and using (\ref{eq:funny-jump-bound})
  concludes the proof.
\end{Proof}

\begin{Cor}[A posteriori bounds for the nonvariational problem including data approximation]
  \label{cor:data-app}
  Let the conditions of Theorem \ref{the:primal-upper} hold. Let $u_h$ now be the the solution of
  (\ref{eq:non-div-fem}) with the following approximations of the bilinear and linear forms
  \begin{gather}
    \bih{u_h}{v_h} := \sum_{K\in\T{}} Q_K^{2k-2} \qp{\frob{\A}{\H(u_h)} v_h} + \sum_{e\in\E} Q_e^{2k} \qp{\sigma h_e \jump{\nabla u_h} \jump{\nabla v_h} + \sigma h_e^{-1} \jump{u_h}\cdot \jump{v_h}}
    \\
    l(v_h) := \sum_{K\in\T{}} Q_K^{2k-2} \qp{f v_h},
  \end{gather}
  then
  \begin{equation}
    \eenorm{u - u_h}{2} 
    \leq 
    C 
    \qp{\sum_{K\in\T{}} \eta_R^2 + \widetilde{\eta_I}^2 + \sum_{e\in\E} \eta_J^2}^{1/2},
  \end{equation}
  and for $\sigma$ large enough
  \begin{equation}
    \sum_{K\in\T{}} \eta_R + \sum_{e\in\E} \eta_J
    \leq
    C \Norm{\Hess u - \Hess_h u_h}_{\leb{2(\W)}} + \sum_{K\in\T{}} \eta_I
  \end{equation} 
  where
  \begin{gather}
    \eta_R^2 := \Norm{P_{k-2} f - P_{k-2}\qp{\frob{\A}{\H(u_h)}}}_{\leb{2}(K)}^2
    \\
    \eta_I^2 := \Norm{\qp{P_{k-2} - \Id}\qp{\frob{\A}{\H(u_h)}}}_{\leb{2}(K)}^2 + \Norm{f - P_{k-2}f}_{\leb{2}(K)}^2
    \\
    \widetilde{\eta_I}^2 := \Norm{\qp{P_{k-2} - \Id}\qp{\frob{\A}{\H(E^2(u_h))}}}_{\leb{2}(K)}^2 + \Norm{f - P_{k-2}f}_{\leb{2}(K)}^2
    \\
    \eta_J^2 := 
    h_e^{-1} \Norm{\jump{\nabla u_h}}_{\leb{2}(e)}^2
    +
    h_e^{-3} \Norm{\jump{u_h}}_{\leb{2}(e)}^2.
  \end{gather}
\end{Cor}
\begin{Proof}
  The proof is the same as that of Theorem \ref{the:primal-upper} with
  the exception that the inconsistency term is slightly more
  complicated. We have that 
  \begin{equation}
    \label{eq:cor-data-app-1}
    \begin{split}
      \E_3 
      &=
      \sup_{v\in \leb{2}(\W), \Norm{v}_{\leb{2}(\W)}\leq 1}
      \frac{1}{\gamma}
      \int_\W P_{k-2}\qp{\frob{\A}{\H(E^2(u_h))}} v - \frob{\A}{\Hess E^2(u_h)} v
      \\
      &=
      \sup_{v\in \leb{2}(\W), \Norm{v}_{\leb{2}(\W)}\leq 1}
      \frac{1}{\gamma}
      \int_\W \bigg[ \qp{P_{k-2}\qp{\frob{\A}{\H(E^2(u_h))}} v - \frob{\A}{\H(E^2(u_h))} v}
      \\
      &\qquad \qquad \qquad \qquad \qquad \qquad \qquad \qquad \qquad 
      + \qp{\frob{\A}{\H(E^2(u_h))} v - \frob{\A}{\Hess E^2(u_h)} v }\bigg].
    \end{split}
  \end{equation}
  The second term in (\ref{eq:cor-data-app-1}) can be bounded using
  the same arguments as in the proof of Theorem
  \ref{the:primal-upper}. The other term is controllable by $\eta_I$
  using Cauchy-Schwarz, concluding the proof.
\end{Proof}

\begin{Rem}[Comparison of the variational and nonvariational estimates]
  Note the different structures of the data approximation terms
  appearing in Lemma \ref{lem:apost-quad} and Corollary
  \ref{cor:data-app}. They are different orders of approximation even
  though the underlying quadrature approximation is the same. This is
  because the estimators control the error in different norms due to
  the different stability frameworks of the underlying PDEs. The data
  approximation terms are of the same order of accuracy as the other
  residual terms using a quadrature of degree $2k-2$. Higher order
  quadrature choice naturally results in the data approximation terms
  becoming higher order, assuming smooth data.

  For the nonvariational scheme the inconsistency is in two distinct
  parts. The first being the nature of the scheme, which was tackled
  in Theorem \ref{the:primal-upper} and is inherently tied to the
  nonconforming nature of the scheme. The second is the quadrature
  approximations. This inconsistency cannot be linked to the
  nonconformity of the scheme, thus another quantity enters into the
  estimator.
\end{Rem}

\begin{Cor}[$\sobh{2}$ a posteriori bounds for the Laplacian]
  Let $\A = -\vec I$ then $\frob{\A}{\Hess u} = -\Delta u$ then the
  following a posteriori bound holds for $u_h$ the approximation given
  by the IP method under quadrature approximation
  (\ref{eq:bilinear-dg-quad}) 
  \begin{equation}
    \eenorm{u - u_h}{2} 
    \leq 
    C 
    \qp{\sum_{K\in\T{}}\eta_R^2 +  \eta_I^2 + \sum_{e\in\E} \eta_J^2}^{1/2},
  \end{equation}
  where
  \begin{gather}
    \eta_R^2 := \Norm{P_{k-2} f + \Delta u_h}_{\leb{2}(K)}^2
    \\
    \eta_I^2 := \Norm{f - P_{k-2} f}_{\leb{2}(K)}^2
    \\
    \eta_J^2 := 
    h_e^{-1} \Norm{\jump{\nabla u_h}}_{\leb{2}(e)}^2
    +
    h_e^{-3} \Norm{\jump{u_h}}_{\leb{2}(e)}^2.
  \end{gather}
  and 
  \begin{equation}
    \sum_{K\in\T{}} \eta_R + \sum_{e\in\E} \eta_J
    \leq 
    C \qp{\Norm{\Hess u - \Hess_h u_h}_{\leb{2}(\W)} + \Norm{f - P_{k-2} f}_{\leb{2}(\W)}}.
  \end{equation}
\end{Cor}
\begin{Proof}
  The Proof consists of observing that the finite element method
  (\ref{eq:non-div-fem}) coincides with the IP method
  (\ref{eq:bilinear-dg-quad}) (albeit with an additional penalty term
  for the gradients) whenever $\A$ is ($\T{}$--wise) constant and
  applying Theorems \ref{the:primal-upper} and \ref{the:primal-lower}.
\end{Proof}

\begin{Rem}[Dual a posteriori error control and the formal dual of (\ref{eq:nondiv})]
  \label{the:dual-upper}
   The formal dual problem is 
  \begin{equation}
    \bid{u}{v} = \int_\W u \frob{\Hess}{\qp{v \A}}
  \end{equation}
  and, to the authors knowledge, there is no general existence result
  for this class of problem. Under the stringent assumption that the
  dual problem is well posed, if $u_h$ solves the nonvariational
  Galerkin method and $u$ is a strong solution then
  \begin{equation}
    \Norm{u - u_h}_{\leb 2(\W)} 
    \leq 
    C 
    \qp{\sum_{K\in\T{}} \eta_{2,R}^2 + \sum_{e \in \E} \eta_{2,J}^2}^{1/2},
  \end{equation}
  where
  \begin{gather*}
    \eta_{2,R}^2 := h_K^4 \Norm{f - \frob{\A}{\H(u_h)}}_{\leb{2}(K)}^2
    \\
    \eta_{2,J}^2 := 
    h_e^{3} \Norm{\jump{\nabla u_h}}_{\leb{2}(e)}^2
    +
    h_e \Norm{\jump{u_h}}_{\leb{2}(e)}^2
  \end{gather*}
  This result is included as, from a
  practical point of view, the estimators seem to be reliable and
  efficient.
\end{Rem}

\section{Numerical experiments}
\label{eq:numerics}

In this section we present some numerical results testing both the
asymptotic behaviour of the estimators derived in previous sections as well as its effectiveness as an indicator for adaptivity.

\subsection{Inconsistency terms from \S\ref{sec:selfadj}}
\label{sec:incon}

\begin{figure}[h!t]
      \caption[]{
        \label{fig:quadsoln}
        The adaptive approximation to the solution of the diffusion problem with problem data given in \S\ref{sec:incon}.
      }
      \centering
      \includegraphics[scale=\figscale, width=0.9\figwidth]
      {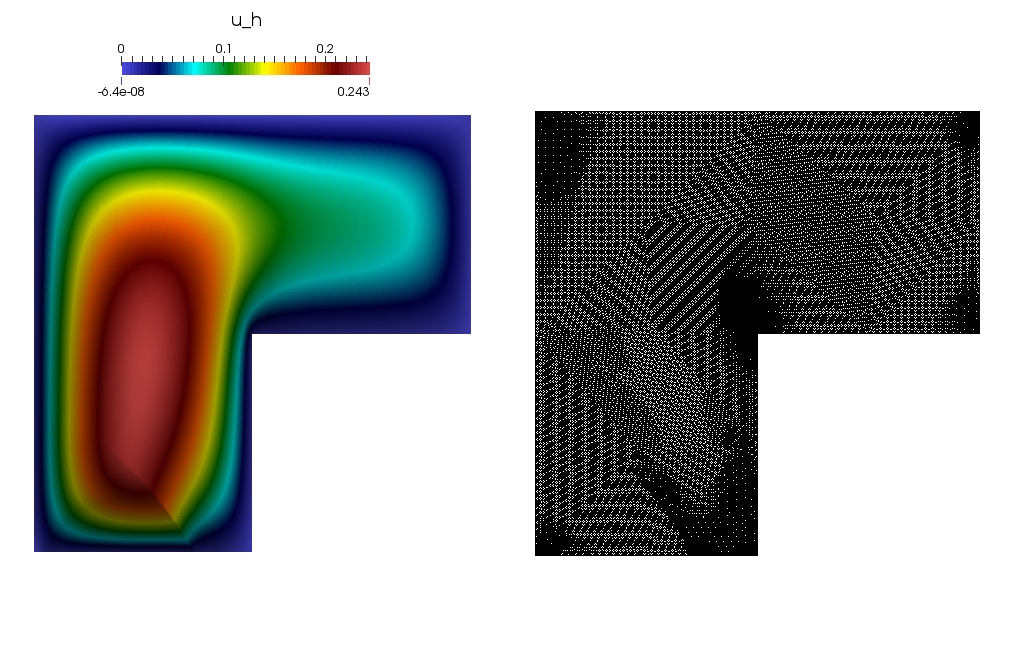} 
    \end{figure}

We begin this section with an adaptive numerical test problem for the
estimator given in Lemma \ref{lem:apost-quad}. We take the diffusion
coefficient to be
\begin{equation}
  \A = 
  \begin{bmatrix}
    1 & 0
    \\
    0 & \arctan{1000\qp{\norm{\qp{x+\frac{1}{2}}^2 + \qp{y-\frac{1}{2}}^2 - 1}}}+2
  \end{bmatrix}
\end{equation}
and $\W = [0,1]^2 - [\frac{1}{2},1]\times [0,\frac{1}{2}]$ to be the standard L-shaped domain.
We use a maximum strategy adaptive algorithm given in Algorithm 1
\cite{alberta} and the estimator given in Lemma
\ref{lem:apost-quad}. We use $k=2$ and take the quadrature degree to
be order $2$. Figure \ref{fig:quadsoln} shows the numerical solution at
the final level of refinement and Figure \ref{fig:quadadapt} shows the
mesh generated and estimator contributions at various levels of
refinement. Note the solution will be singular at the nonconvex corner
of the domain. In addition the problem data is chosen such that $\A$
has a very large gradient over a large portion of the domain which
induces a further singularity. Both seem to be captured correctly
using this algorithm.
    \begin{algorithm}[ht]
    \caption{Maximum strategy type adaptive algorithm incorporating inconsistencies.}
    \begin{algorithmic}
    \label{Alg}
      \State Given constants $C_1, C_2, C_3$ such that $C_1 + C_2+ C_3 = 1$, \tol and $\threshhold, \coarsethresh$;
      \State Set N = 0;
      \While{$C_1\Eta_R + C_2\Eta_I + C_3\Eta_J > \tol$ AND N $\leq$ Max Iterations}
      \ForAll{$K\in\T{}$}
      \State Compute $\eta_{R}, \eta_{J}$ and $\eta_{I}$;
      \State $\Eta_R = \Eta_R + \eta_{R}^2$;
      \State $\Eta_I = \Eta_I + \eta_{I}^2$;
      \ForAll{$e\in K$}
      \State $\Eta_J = \Eta_J + \frac{1}{2}\eta_{J,e}^2$;
      \EndFor
      \State Set $\eta_K = \eta_{R} + \frac{1}{2}\eta_{J} + \eta_{I}$;
      \If{$ \eta_K \geq \threshhold \max_{L\in\T{}}{(\eta_{L})}$}
      \State Mark element $K$ for refinement;
      \EndIf
      \If{$\eta_{K} \leq \coarsethresh \max_{L\in\T{}}{(\eta_{L})}$}
      \State Mark element $K$ for coarsening;
      \EndIf
      \EndFor
      \State $\Eta_R = \sqrt{\Eta_R}$, $\Eta_J = \sqrt{\Eta_J}$, $\Eta_I = \sqrt{\Eta_I}$;
      \State Perform any mesh change;
      \State N = N+1;
      \EndWhile
   \end{algorithmic}
    \end{algorithm}

\begin{figure}[h!t]
  \caption[]
  {\label{fig:quadadapt} Testing the estimator given in Lemma
    \ref{lem:apost-quad} as a driver for adaptivity. The problem data
    is described in \S \ref{sec:incon}. We consider various
    iterates of the adaptive procedure looking at the mesh and the
    estimator components. The top estimator is the inconsistency terms
    and the bottom the interior and jump residual. Notice the
    inconsistency term is extremely localised to where the diffusion
    coefficient has large gradient and becomes comparatively
    negligible once the mesh is sufficiently resolved in that region.}
  \subfigure[][Iterate 3] {
    \includegraphics[scale=\figscale, width=0.45\figwidth]
                    {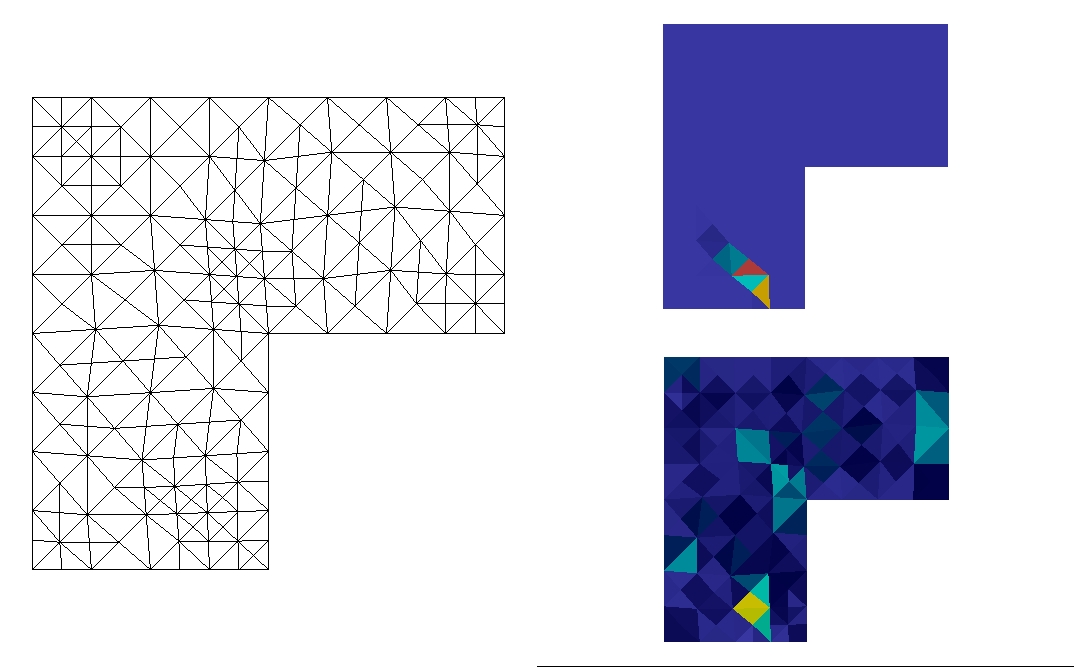} }
  \subfigure[][Iterate 4] {
    \includegraphics[scale=\figscale, width=0.45\figwidth]
    {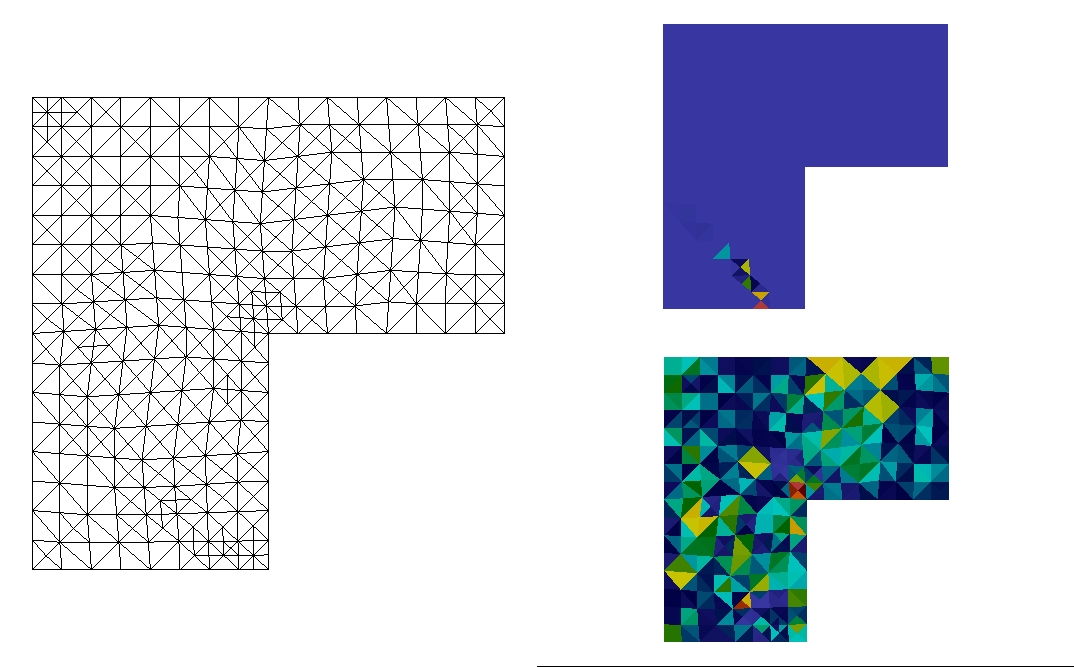}                  }             
  \subfigure[][Iterate 5] {
    \includegraphics[scale=\figscale, width=0.45\figwidth]
                    {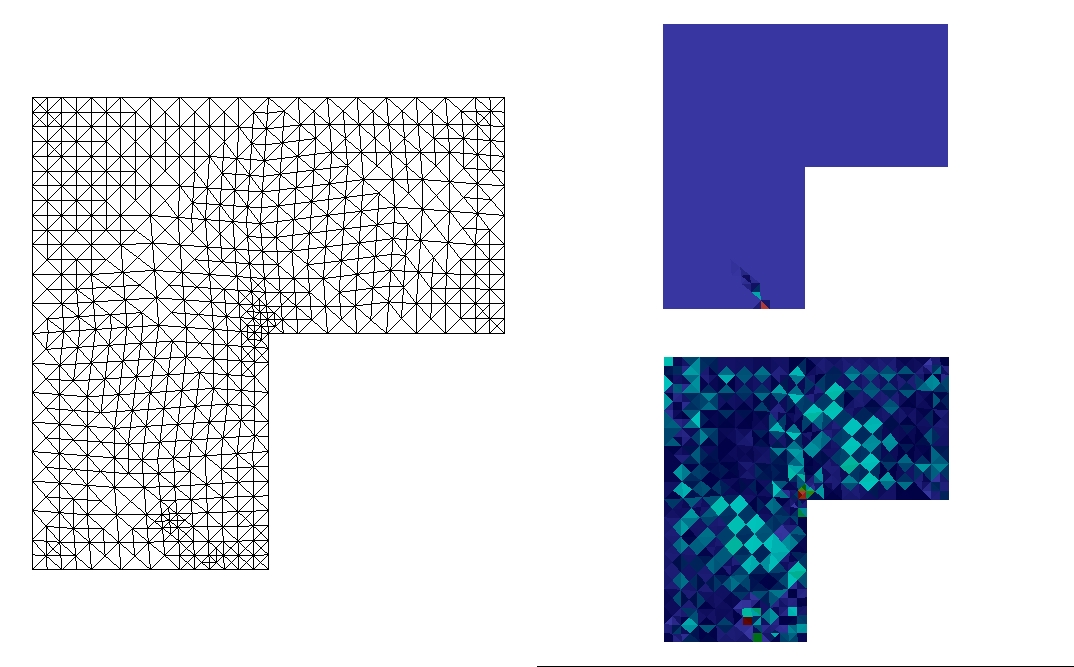} }
  \subfigure[][Iterate 6] {
    \includegraphics[scale=\figscale, width=0.45\figwidth]
    {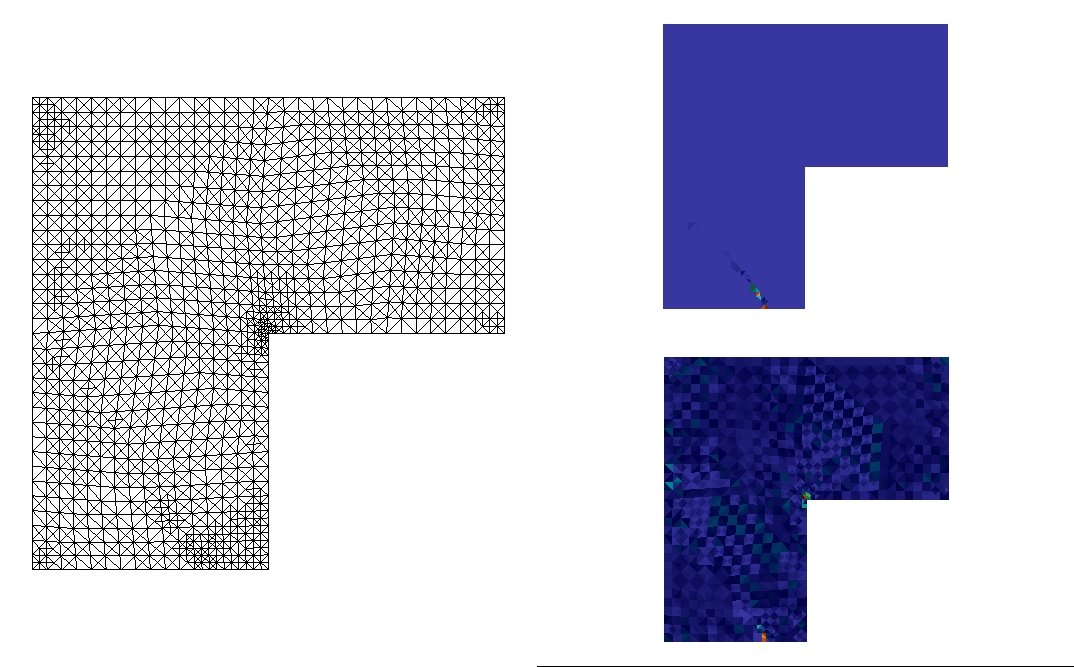}                   }            
  \subfigure[][Iterate 7] {
    \includegraphics[scale=\figscale, width=0.45\figwidth]
                    {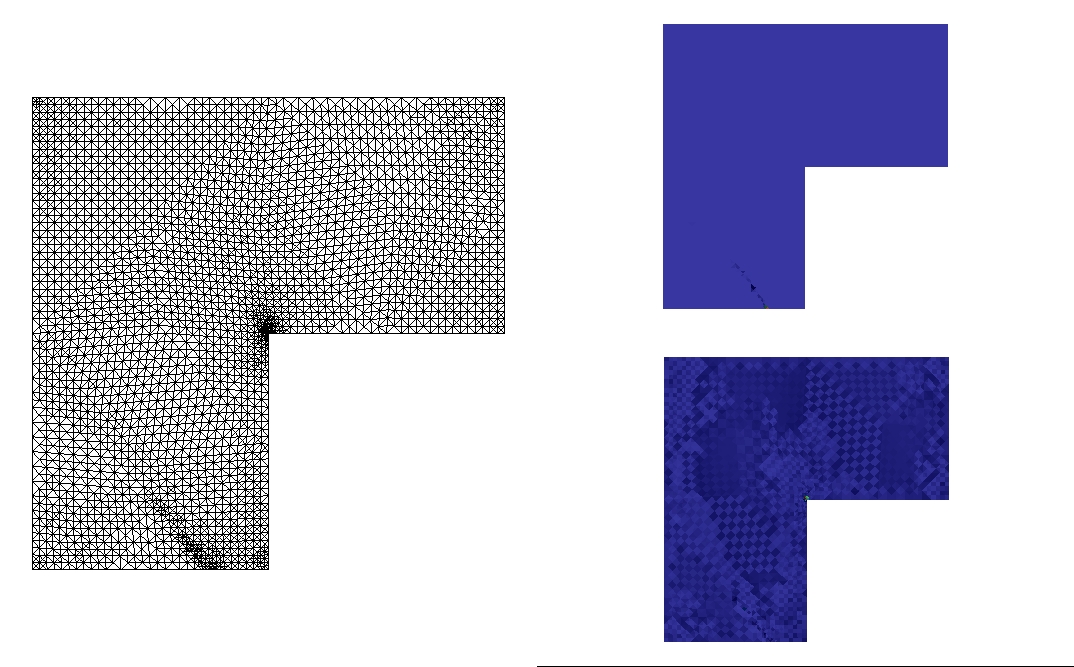} }
  \subfigure[][Iterate 8] {
    \includegraphics[scale=\figscale, width=0.45\figwidth]
    {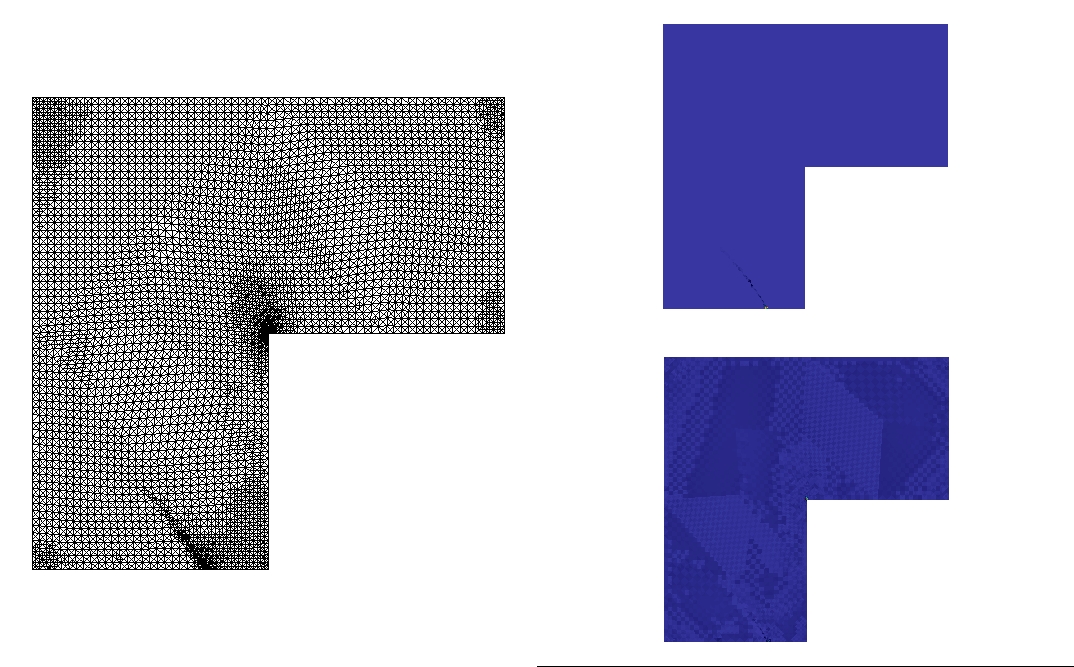}                    }           
  
\end{figure}

\clearpage

\subsection{Inconsistent nonvariational approximations from \S\ref{sec:unbalanced}}
\label{sec:nonadapt}

We are studying both the primal and dual
estimators derived in \S\ref{sec:unbalanced} for nonconforming,
inconsistent schemes approximating nonvariational problems. We examine the asymptotic behaviour of
the estimators from Theorems \ref{the:primal-upper} and
\ref{the:dual-upper} and also their effectiveness as drivers for adaptive algorithms.

\subsubsection{Test 1 : Asymptotic behaviour} 
\begin{figure}[h!]
  \caption[]
  {\label{fig:nonadapt} Testing the asymptotic behaviour of the
    estimator given in Theorem \ref{the:primal-upper}. The problem
    data is described in \S \ref{sec:nonadapt} - Test 1.}
  \subfigure[][When $u\in\cont{\infty}(\W)$ is given by
  (\ref{eq:u-smooth})] {
    \includegraphics[scale=\figscale, width=0.47\figwidth]
    {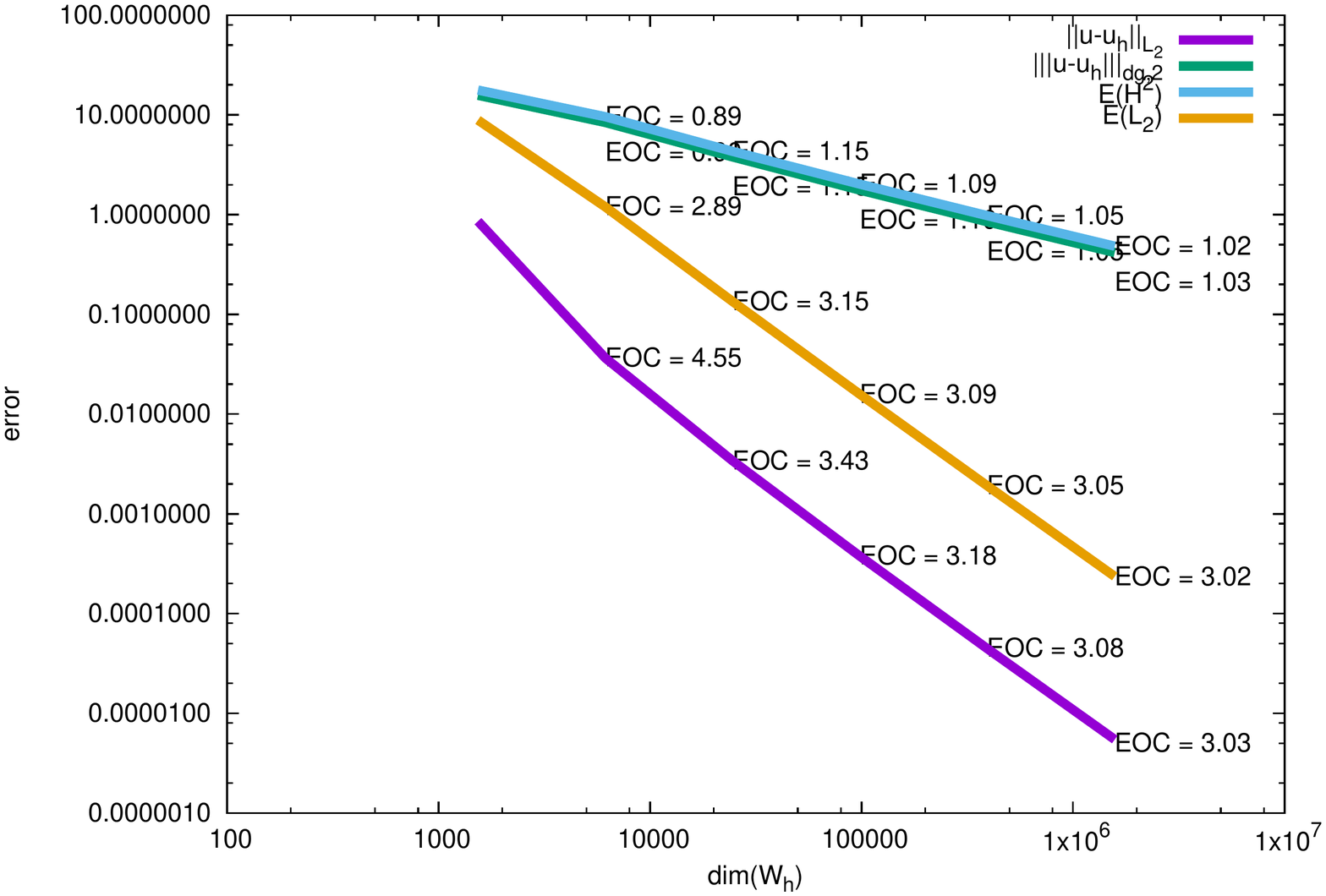} } \hfill \subfigure[][When
  $u\in\sobh2(\W)\not \ \ \sobh{3}(\W)$ is given by
  (\ref{eq:u-notsosmooth})] {
    \includegraphics[scale=\figscale, width=0.47\figwidth]
    {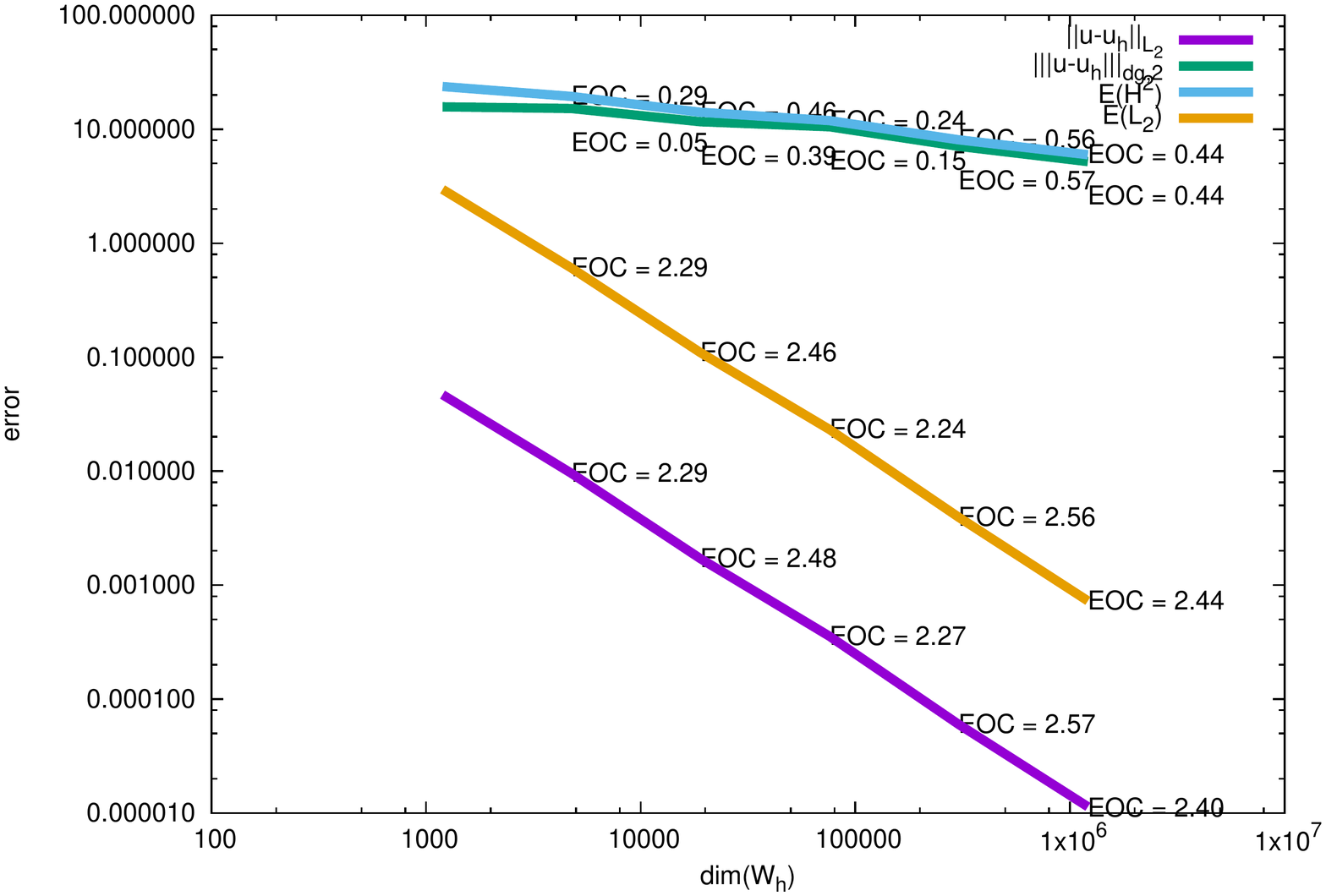}                                
  }
\end{figure}

With $d=2$ if
$\frob{\A}{\Hess u}$ is a uniformly elliptic operator then the
Campanato condition for strong solutions is satisfied
\cite{Tarsia:2004} and Assumption \ref{ass:cont-inf-sup} is
satisfied. We fix $\W = \qb{\frac{1}{2},\frac{1}{2}}^2$ use the
problem data
\begin{equation}
  \label{eq:diffusion}
  \A = 
  \begin{bmatrix}
    2 & 0
    \\
    0 & \sin{2\pi x}\sin{2\pi y}+2.
  \end{bmatrix}
\end{equation}
We have chosen $\A$ in this fashion such that it varies slowly over
the domain and the inconsistency term arising from quadrature
approximation is negligible. We choose $f$ such that the exact
solution is given by either
\begin{equation}
  \label{eq:u-smooth}
  u(\vec x) = \sin{2\pi x} \sin{2\pi y} \in \cont{\infty}(\W)
\end{equation}
or 
\begin{equation}
  \label{eq:u-notsosmooth}
  u(\vec x) =
  \begin{cases}
    \frac{1}{4}\qp{\cos{8\pi\norm{\vec x - \frac{1}{2}}^2}+1} \text{ if } \norm{\vec x - \frac{1}{2}}^2 \leq \frac{1}{8}
    \\
    0
  \end{cases}
  \in \sobh 2(\W) \not \ \  \sobh3(\W).
\end{equation}
In Figure \ref{fig:nonadapt} we summarise the results of this test and
demonstrate that computationally the asymptotic convergence rate of the estimator is the same
as that of the error as predicted by the theory given in Theorems
\ref{the:primal-upper} and \ref{the:primal-lower}.

\subsubsection{Test 2 : Adaptive algorithms}

\begin{figure}[h!]
  \caption[]
  {\label{fig:nonvaradapt} Testing the estimator given in Theorem
    \ref{the:primal-upper} as a driver for adaptivity. The problem
    data is described in \S\ref{sec:nonadapt} - Test 2 such that the
    solution is given by (\ref{eq:u-notsosmooth}). We consider various
    iterates of the adaptive procedure governed by Algorithm 1 looking
    at the mesh generated. Notice the mesh is refined around where the
    solution is nonsmooth and the convergence rate of the adaptive
    algorithm is quicker than the uniform counterpart from Test 1.}
  \subfigure[][Iterate 1.] {
    \includegraphics[scale=\figscale, width=0.45\figwidth]
    {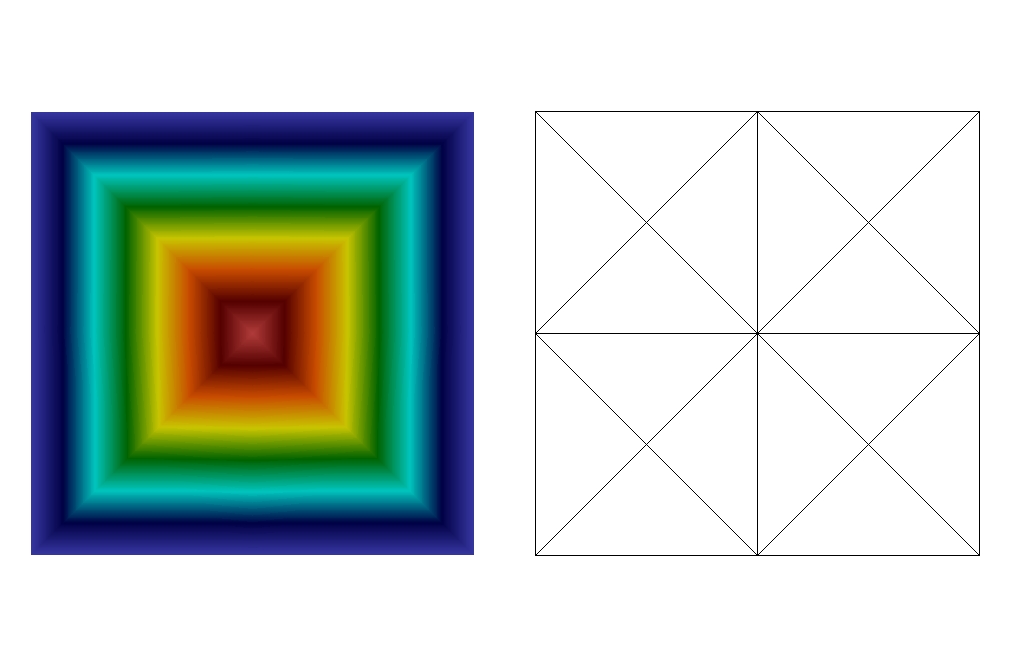} }
  \subfigure[][Iterate 11.] {
    \includegraphics[scale=\figscale, width=0.45\figwidth]
    {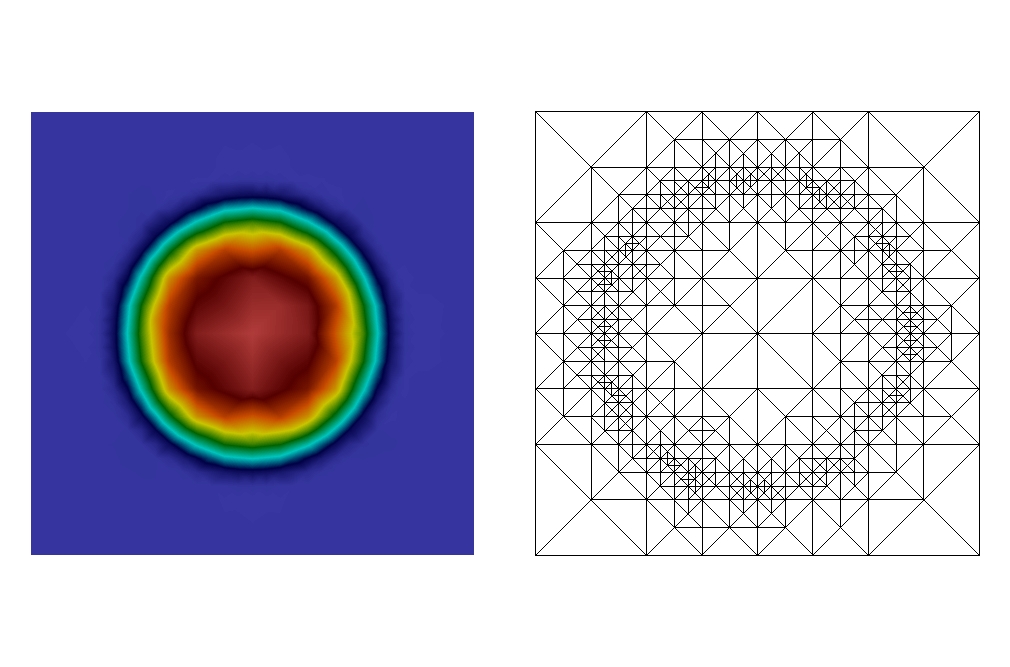}                  }             
  \subfigure[][Iterate 13.] {
    \includegraphics[scale=\figscale, width=0.45\figwidth]
    {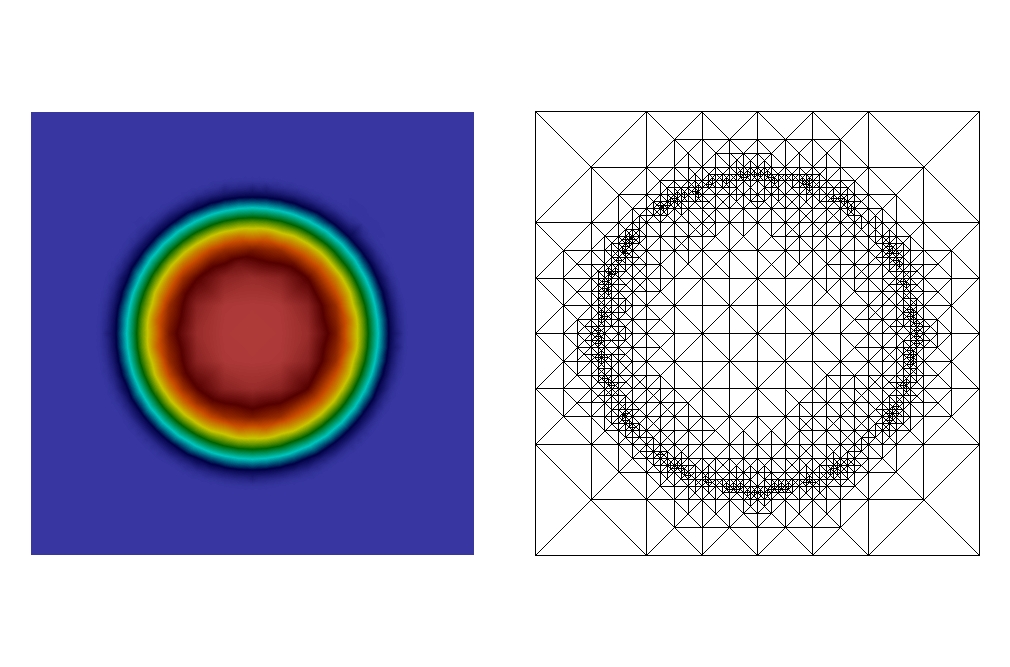} }
  \subfigure[][Iterate 18.] {
    \includegraphics[scale=\figscale, width=0.45\figwidth]
    {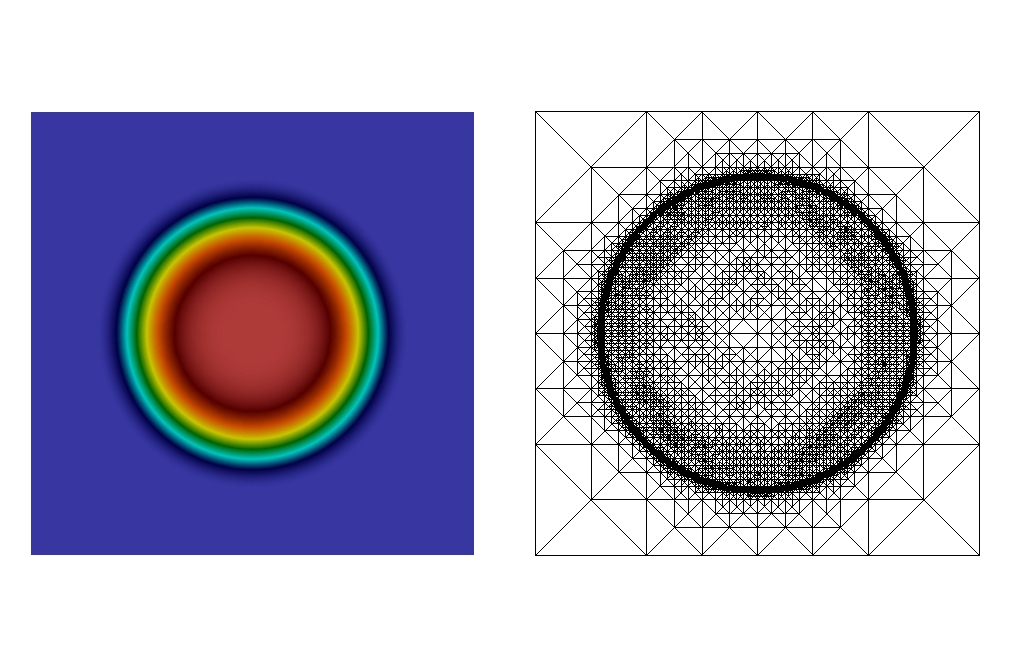}                    }           
  \subfigure[][Convergence rates for the adaptive algorithm.] {
    \includegraphics[scale=\figscale, width=0.55\figwidth]
    {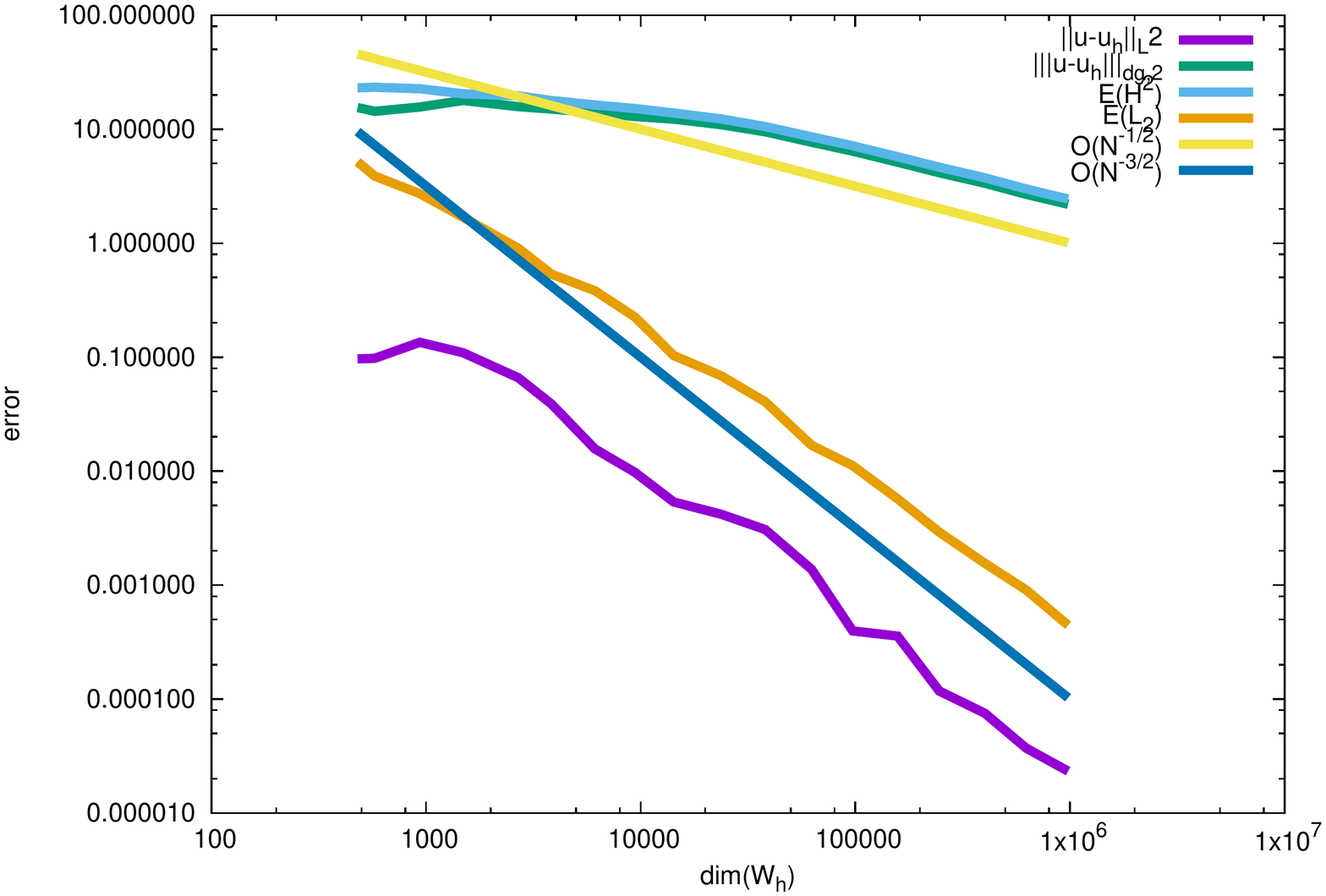}                    }           
  
\end{figure}

We fix $\W = [\frac{1}{2}, \frac{1}{2}]^2$, as in Test 1, and choose
$\A$ as in (\ref{eq:diffusion}). We choose $f$ such that the exact
solution is given in (\ref{eq:u-notsosmooth}). Note this function is
$\sobh{2}(\W) \not \ \ \sobh3(\W)$. We use the adaptive algorithm
given in Algorithm 1 and, in Figure \ref{fig:nonvaradapt}, examine the
convergence of said algorithm. The rates are increased from those
given in Figure \ref{fig:nonadapt} (approximately $\Oh(N^{-1/4})$) to $\Oh(N^{-1/2})$ when measured in
the $\sobh2(\W)$ norm.

\subsubsection{Test 3 : Quadrature effects, smooth solution}

We fix $\W = [\frac{1}{2}, \frac{1}{2}]^2$, as in Test 1 and 2, and
choose 
\begin{equation}
  \A = 
  \begin{bmatrix}
    1 & 0
    \\
    0 & a_{11}(\vec x),
  \end{bmatrix}
\end{equation}
where
\begin{equation}
  a_{11}(\vec x) 
  =
  \begin{cases}
    10\sin{100\pi x}\sin{100\pi y} + 11 \text{ if } x, y \geq 0
    \\
    11 \text{ otherwise}.
  \end{cases}
\end{equation}
Notice the diffusion coefficient oscillates much more heavily in the
upper right quadrant of $\W$ than that of the previous tests. The is
to force the inconsistency term arising from quadrature to have a
larger impact. We choose $f$ such that the exact solution is smooth
and given in (\ref{eq:u-smooth}). We use the adaptive algorithm given
in Algorithm 1 and, in Figure \ref{fig:nonvaradapt}, examine the
behaviour of the inconsistency term compared to the rest of the
estimator. The increase in size and oscillations of the diffusion term
force the inconsistency term in the estimator to dominate at the initial coarse
mesh scale which results in a denser mesh in the upper right quadrant
of $\W$.

\begin{figure}[h!t]
  \caption[]
  {\label{fig:quadadapt} Testing the estimator given in Theorem
    \ref{the:primal-upper} as a driver for adaptivity. The problem
    data is described in \S \ref{sec:nonadapt} - Test 3. In this case the
    solution is smooth and the diffusion coefficient oscillates
    heavily in the upper right quadrant. We consider various iterates
    of the adaptive procedure looking at the mesh and the estimator
    components. The top estimator is the inconsistency terms and the
    bottom the interior and jump residual. Notice after the mesh is
    sufficiently resolved where the oscillations are heavy, the
    inconsistency term becomes comparatively negligible.}
  \subfigure[][Iterate 5] {
    \includegraphics[scale=\figscale, width=0.37\figwidth]
    {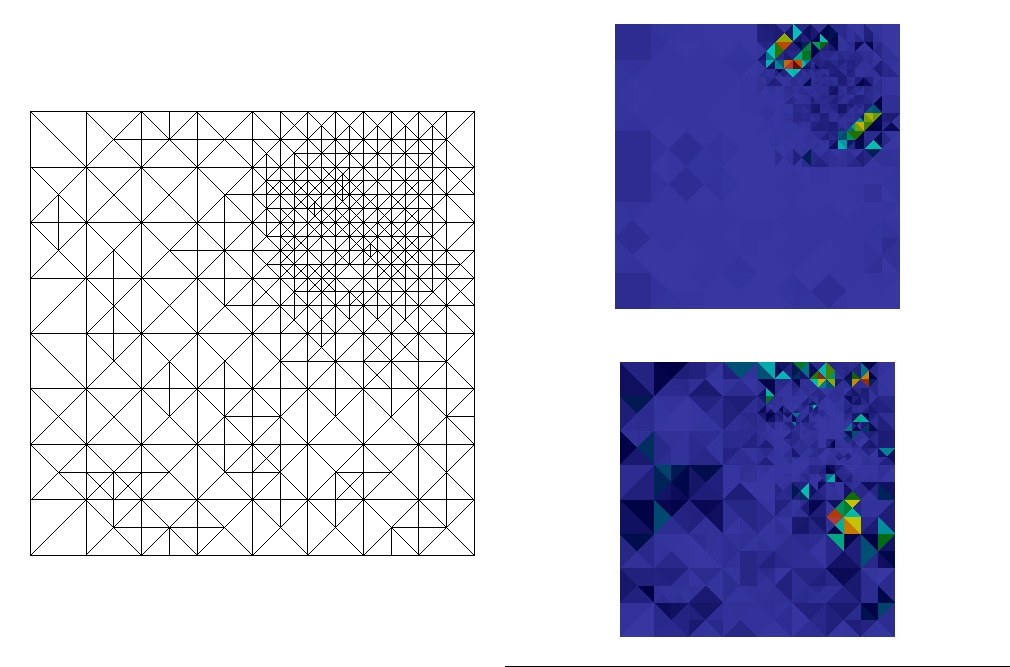} }
  \subfigure[][Iterate 6] {
    \includegraphics[scale=\figscale, width=0.37\figwidth]
    {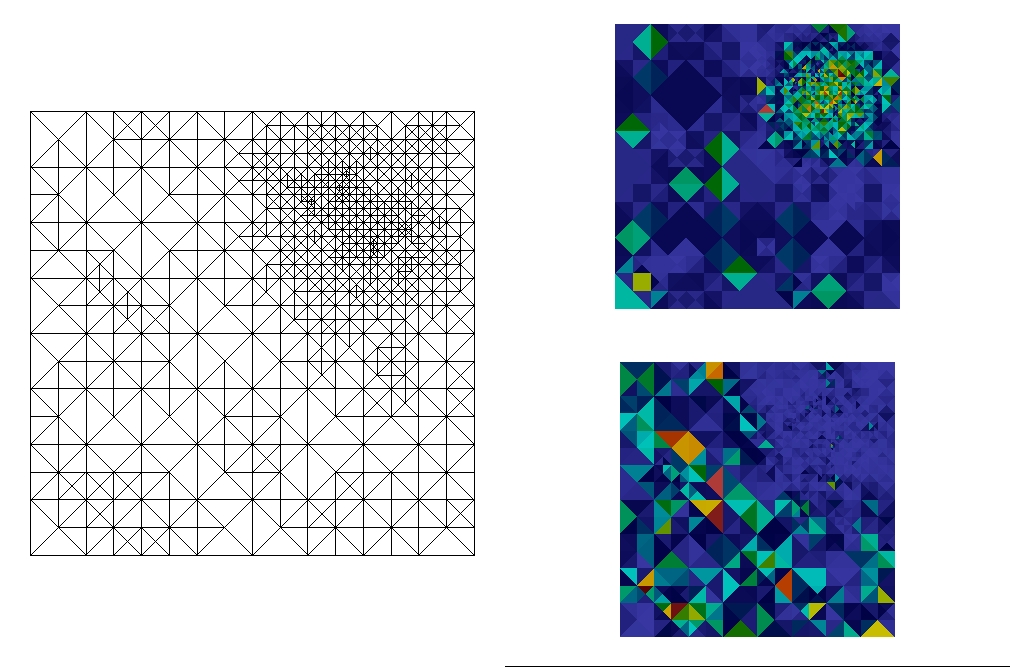}                  }             
  \subfigure[][Iterate 7] {
    \includegraphics[scale=\figscale, width=0.37\figwidth]
    {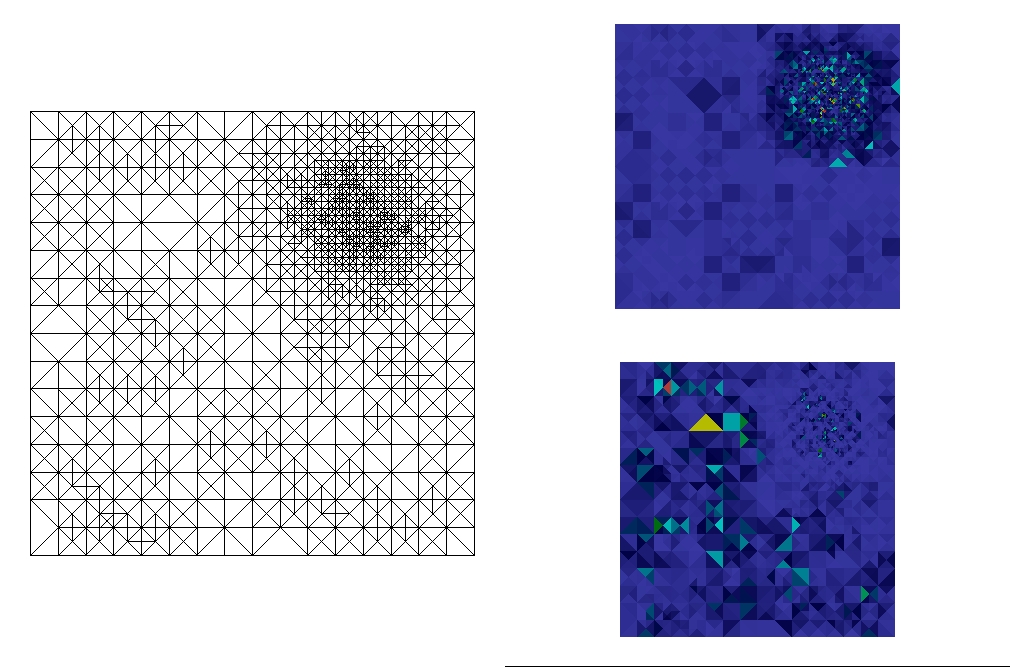} }
  \subfigure[][Iterate 8] {
    \includegraphics[scale=\figscale, width=0.37\figwidth]
    {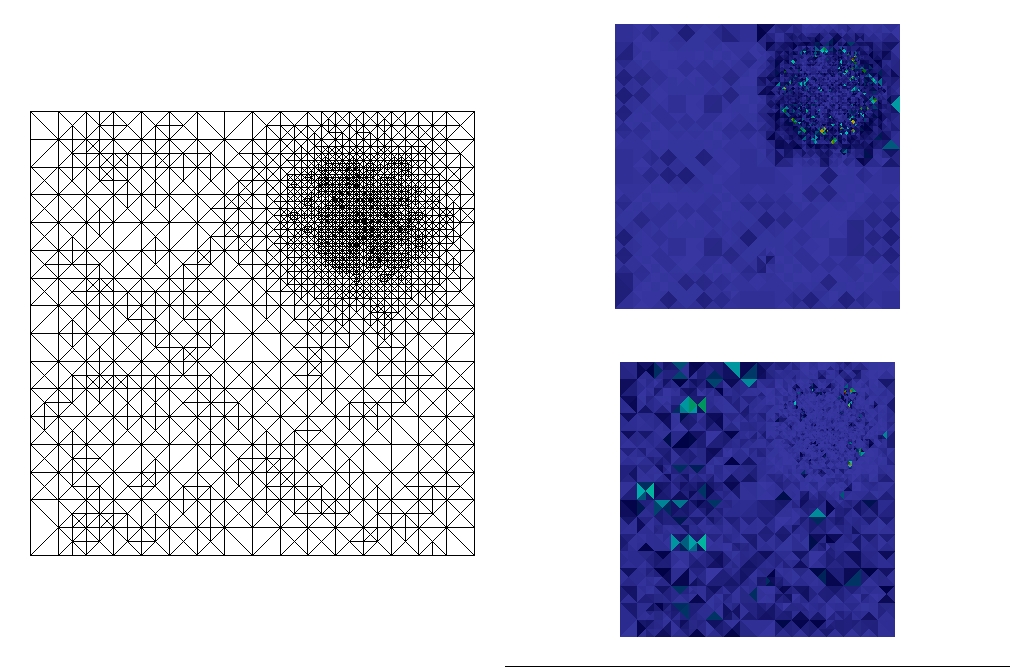}                   }            
  \subfigure[][Iterate 9] {
    \includegraphics[scale=\figscale, width=0.37\figwidth]
    {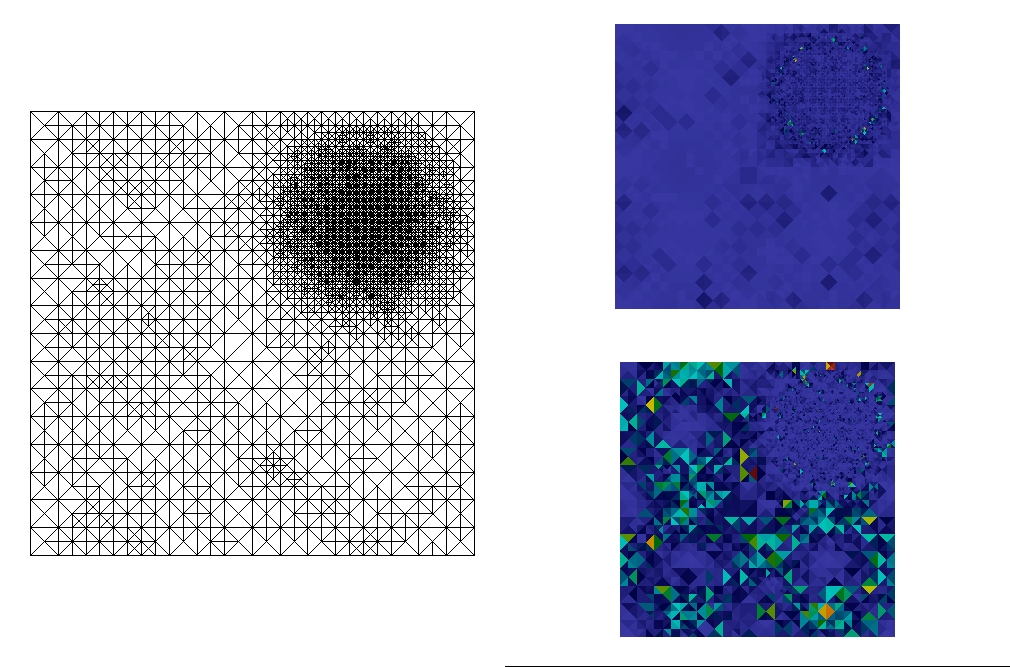} }
  \subfigure[][Iterate 10] {
    \includegraphics[scale=\figscale, width=0.37\figwidth]
    {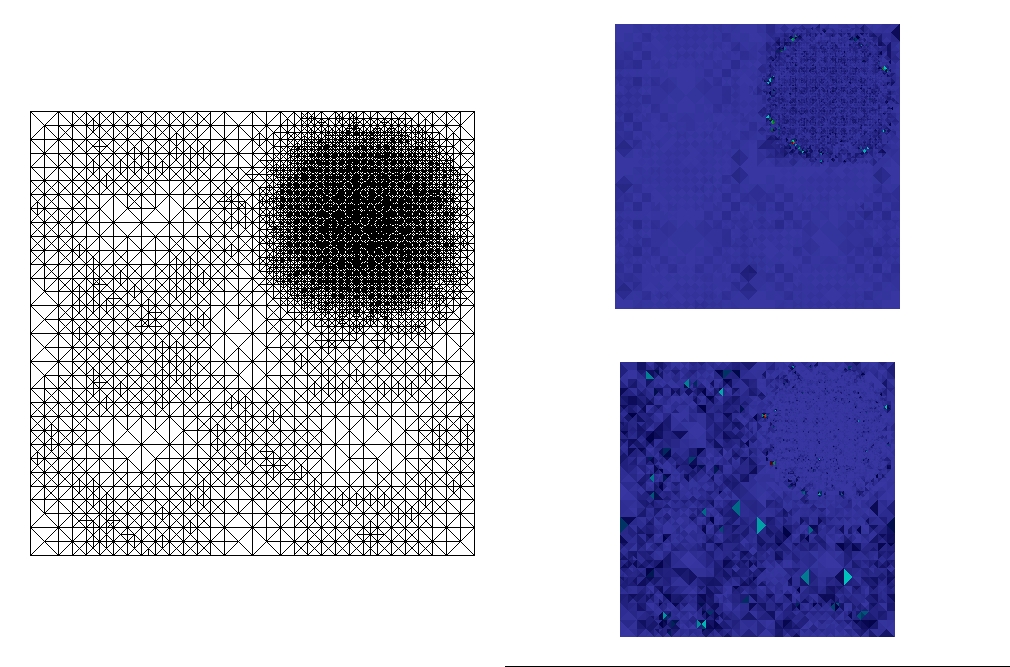}                    }           
  
\end{figure}

\subsubsection{Test 4 : Quadrature effects, nonsmooth solution}

\begin{figure}[h!t]
  \caption[]
  {\label{fig:quadadaptnonsmooth} Testing the estimator given in
    Theorem \ref{the:primal-upper} as a driver for adaptivity. The
    problem data is described in \S \ref{sec:nonadapt} - Test 4. In
    this case the solution is not smooth and the diffusion coefficient
    oscillates heavily in the upper right quadrant. We consider
    various iterates of the adaptive procedure looking at the mesh and
    the estimator components. The top estimator is the inconsistency
    terms and the bottom the interior and jump residual. Notice the
    mesh is very resolved around where the solution is not smooth and also
    where the oscillations in the problem data are heavy.}
  \subfigure[][Iterate 5] {
    \includegraphics[scale=\figscale, width=0.37\figwidth]
    {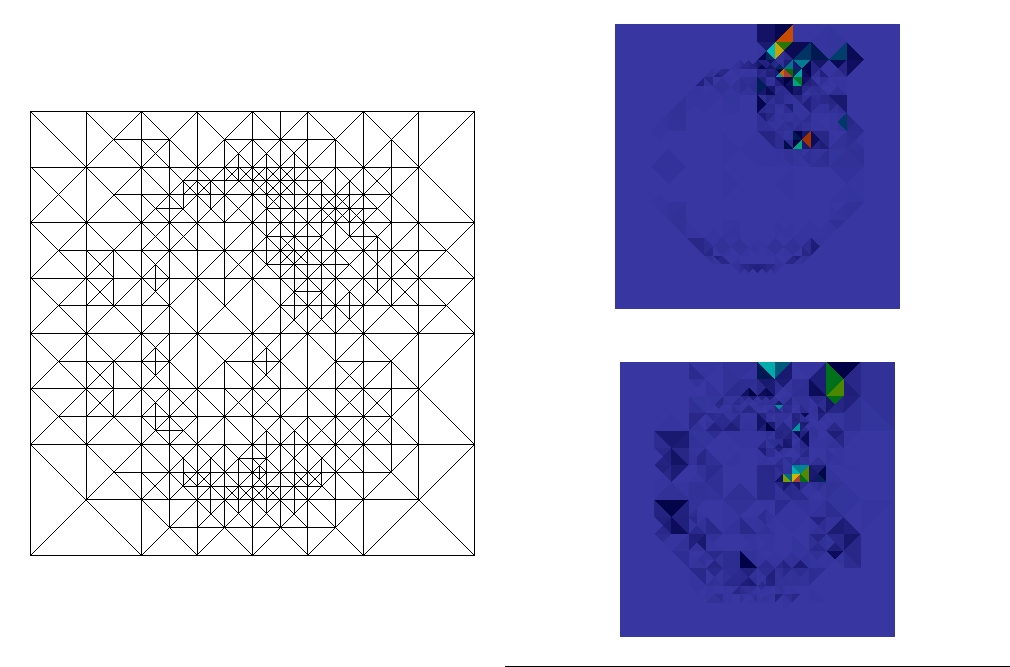} }
  \subfigure[][Iterate 6] {
    \includegraphics[scale=\figscale, width=0.37\figwidth]
    {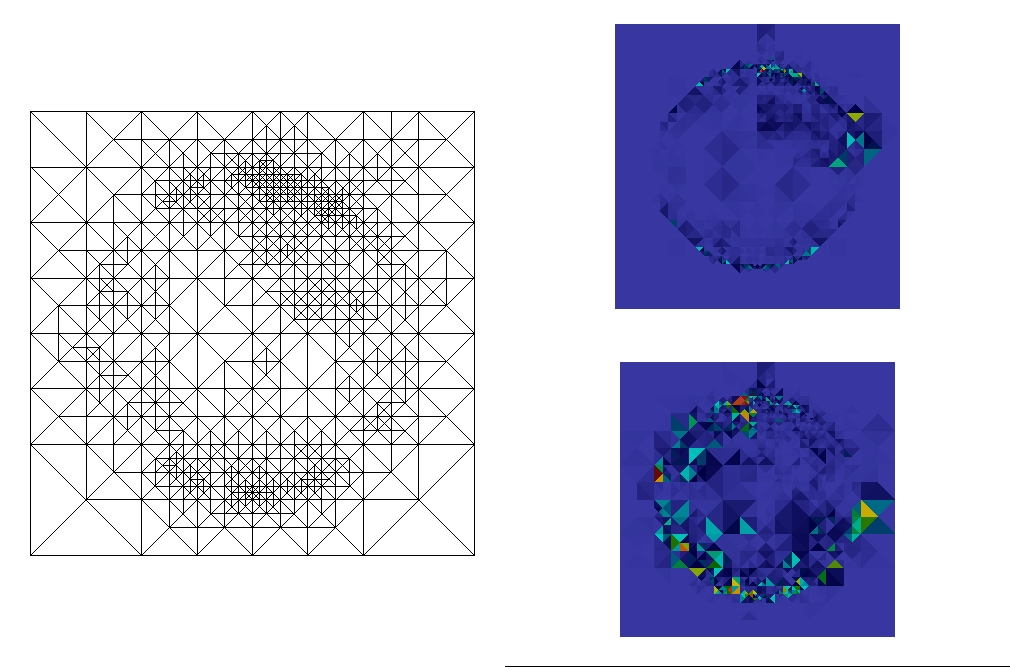}                  }             
  \subfigure[][Iterate 7] {
    \includegraphics[scale=\figscale, width=0.37\figwidth]
    {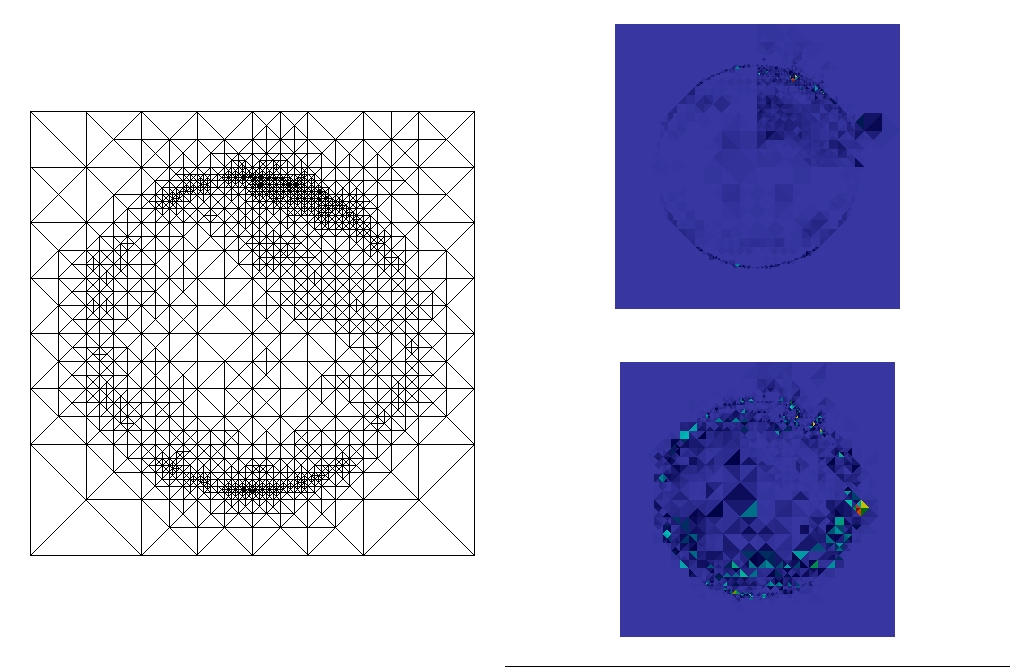} }
  \subfigure[][Iterate 8] {
    \includegraphics[scale=\figscale, width=0.37\figwidth]
    {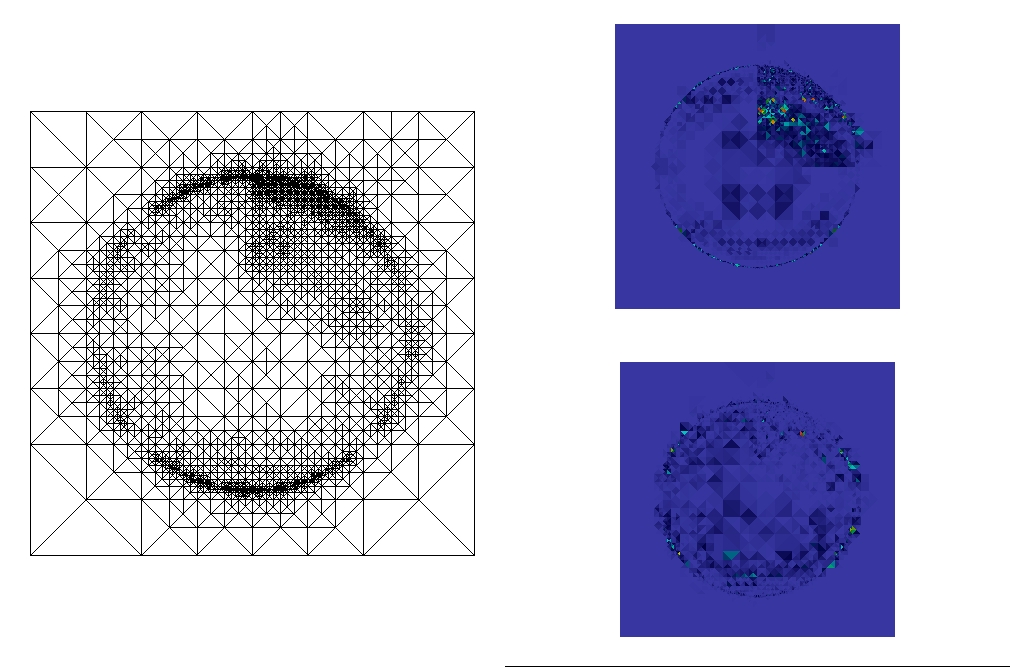}                   }            
  \subfigure[][Iterate 9] {
    \includegraphics[scale=\figscale, width=0.37\figwidth]
    {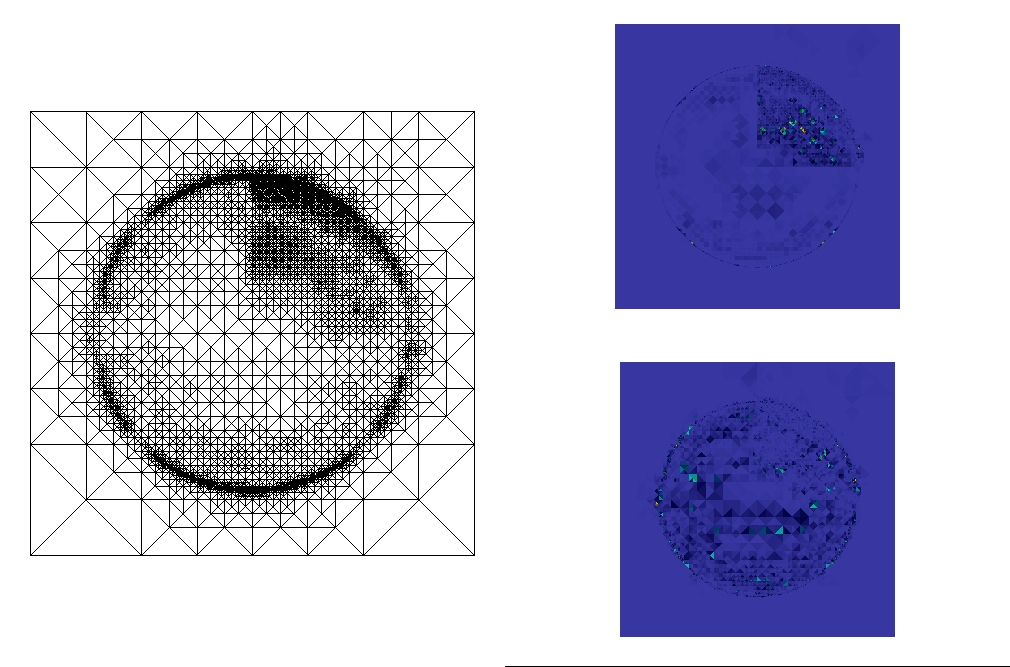} }
  \subfigure[][Iterate 10] {
    \includegraphics[scale=\figscale, width=0.37\figwidth]
    {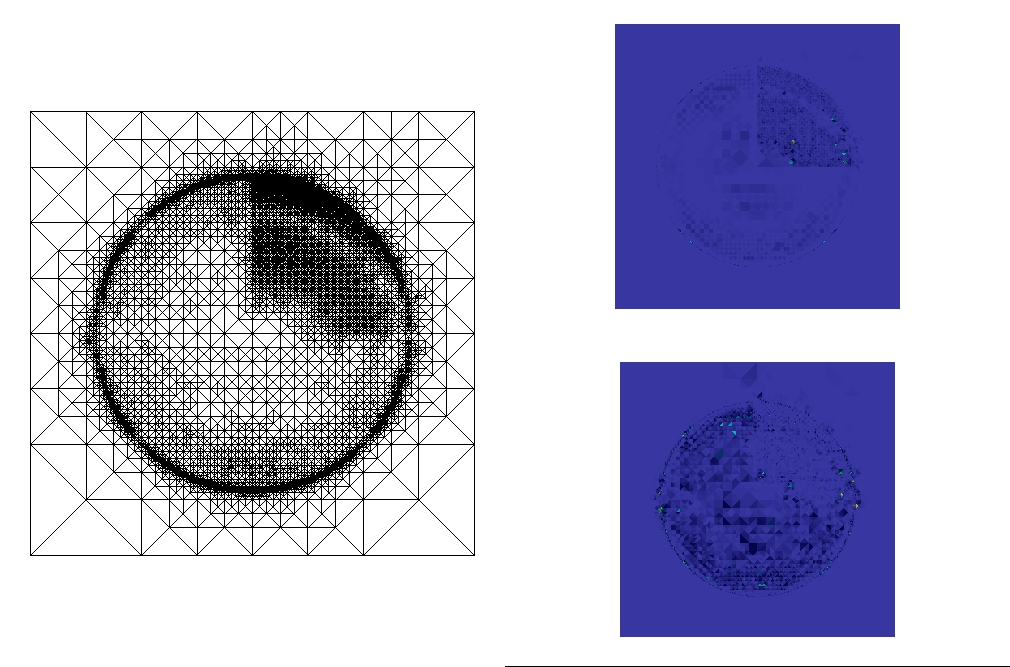}                    }           
  
\end{figure}

This test is the same as Test 4, with the exception that the solution
is given by (\ref{eq:u-notsosmooth}). Figure
\ref{fig:quadadaptnonsmooth} summaries the results.

\bibliographystyle{alpha}
\bibliography{./tristansbib,./tristanswritings}

\end{document}